\documentclass{amsart}

\usepackage{amsfonts,amsmath,amssymb,verbatim}

\theoremstyle{plain}
\newtheorem{theorem}{Theorem}[section]
\newtheorem{lemma}[theorem]{Lemma}
\newtheorem{proposition}[theorem]{Proposition}
\newtheorem{corollary}[theorem]{Corollary}

\theoremstyle{definition}
\newtheorem{definition}[theorem]{Definition}

\numberwithin{equation}{section}

\newcommand{\qbin}[2]{\genfrac{[}{]}{0pt}{}{#1}{#2}}
\newcommand{\qbins}[2]{{\textstyle\genfrac{[}{]}{0pt}{}{#1}{#2}}}
\newcommand{\Z}{\mathbb{Z}}
\newcommand{\Zp}{\mathbb{Z}_+}
\newcommand{\N}{\mathbb{N}}
\newcommand{\Q}{\mathbb{Q}}

\newcommand{\I}{\mathcal{I}}
\newcommand{\C}{\mathcal{C}}
\newcommand{\B}{\mathcal{B}}
\renewcommand{\P}{\mathcal{P}}
\newcommand{\bs}{\boldsymbol}
\newcommand{\vm}{\bs m}
\newcommand{\vn}{\bs n}
\newcommand{\vQ}{\bs Q}
\newcommand{\vA}{\bs A}
\newcommand{\vu}{\bs u}
\newcommand{\vv}{\bs v}
\newcommand{\ve}{\bs e}
\newcommand{\y}{\overline{y}}
\newcommand{\lb}{\overline{l}}

\begin{document}

\title[Conjugate Bailey pairs]{Conjugate Bailey pairs. \\[1mm]
{\tiny From configuration sums and fractional-level string 
functions to Bailey's lemma}}
 
\author[A.~Schilling]{Anne Schilling}
\address{Instituut voor Theoretische Fysica, Universiteit van Amsterdam,
Valckenierstraat 65, 1018 XE Amsterdam, The Netherlands}
\email{schillin@wins.uva.nl, warnaar@wins.uva.nl}
\author[S.~O.~Warnaar]{S.~Ole Warnaar}
 
\subjclass{Primary 05A30, 05A19; Secondary 82B23, 17B67, 33D90}
 
\begin{abstract}
In this paper it is shown that the one-dimensional configuration sums
of the solvable lattice models of Andrews, Baxter and Forrester and the
string functions associated with admissible representations of the
affine Lie algebra A$_1^{(1)}$ as introduced by Kac and Wakimoto can
be exploited to yield a very general class of conjugate Bailey pairs.
Using the recently established fermionic or constant-sign expressions
for the one-dimensional configuration sums, our result is employed to derive
fermionic expressions for fractional-level string functions, 
parafermion characters and A$_1^{(1)}$ branching functions.
In addition, $q$-series identities are obtained whose Lie algebraic
and/or combinatorial interpretation is still lacking.
\end{abstract}
 
\maketitle

\setcounter{section}{-1}

\section{Notation}
Throughout the paper the following notation is used.
$\N$ are the positive integers, $\Zp$ the nonnegative integers,
$\N_p=\{1,\dots,p\}$, $\Z_p=\{0,\dots,p-1\}$.
For $n\in\Z$, $\binom{n}{2}=n(n-1)/2$.

\section{The Bailey lemma}\label{sec BL}

In an attempt to clarify Rogers' second proof \cite{Rogers17} of the 
Rogers--Ramanujan identities, Bailey \cite{Bailey49} was led to the 
following simple but important observation.
\begin{lemma}
If $\alpha=\{\alpha_L\}_{L\geq 0},\dots,\delta=\{\delta_L\}_{L\geq 0}$
are sequences that satisfy
\begin{equation}\label{Bp}
\beta_L=\sum_{r=0}^L \alpha_r u_{L-r}v_{L+r} \quad \text{and} \quad
\gamma_L=\sum_{r=L}^{\infty} \delta_r u_{r-L} v_{r+L},
\end{equation}
then
\begin{equation}\label{abcd}
\sum_{L=0}^{\infty} \alpha_L \gamma_L
=\sum_{L=0}^{\infty} \beta_L \delta_L.
\end{equation}
\end{lemma}
The proof is straightforward and merely requires an interchange of sums.
Of course, in the above suitable convergence conditions need to be imposed
to make the definition of $\gamma$ and the interchange of sums meaningful.

In applications of his transform, Bailey chose $u_L=1/(q)_L$ and
$v_L=1/(aq)_L$, with the usual definition of the $q$-raising factorial,
\begin{equation*}
(a)_{\infty}=(a;q)_{\infty}=\prod_{k=0}^{\infty}(1-aq^k)
\end{equation*}
and
\begin{equation*}
(a)_L=(a;q)_L=(a)_{\infty}/(aq^L)_{\infty}
\end{equation*}
for all $L\in\Z$.
With this choice, equation \eqref{Bp} reads
\begin{equation}\label{BP}
\beta_L=\sum_{r=0}^L \frac{\alpha_r}{(q)_{L-r}(aq)_{L+r}}
\end{equation}
and
\begin{equation}\label{CBP}
\gamma_L=\sum_{r=L}^{\infty} \frac{\delta_r}{(q)_{r-L}(aq)_{r+L}}.
\end{equation}
A pair of sequences that satisfies \eqref{BP} is called a Bailey pair
relative to $a$. Similarly, a pair satisfying \eqref{CBP} is called
a conjugate Bailey pair relative to $a$.

Still following Bailey, one can employ the $q$-Saalsch\"utz summation
\cite[Eq. (II.12)]{GR90} to establish that $(\gamma,\delta)$ with
\begin{equation}\label{gd}
\begin{split}
\gamma_L&=\frac{(\rho_1)_L(\rho_2)_L(aq/\rho_1\rho_2)^L}
{(aq/\rho_1)_L(aq/\rho_2)_L}
\frac{1}{(q)_{M-L}(aq)_{M+L}} \\[2mm]
\delta_L&=\frac{(\rho_1)_L(\rho_2)_L(aq/\rho_1\rho_2)^L}
{(aq/\rho_1)_M(aq/\rho_2)_M}
\frac{(aq/\rho_1\rho_2)_{M-L}}{(q)_{M-L}}
\end{split}
\end{equation}
provides a conjugate Bailey pair.

Unfortunately, Bailey outrightly rejected the above conjugate Bailey pair 
as too complicated to yield any results of interest and focussed on the
simpler case obtained by letting $M$ go to infinity.
Doing so as well as letting the indeterminates $\rho_1$ and $\rho_2$
tend to infinity yields
\begin{equation}\label{gdinf}
\gamma_L=\frac{a^L q^{L^2}}{(aq)_{\infty}} \quad \text{and} \quad
\delta_L=a^L q^{L^2},
\end{equation}
which substituted into \eqref{abcd} gives
\begin{equation}\label{ab}
\frac{1}{(aq)_{\infty}}\sum_{L=0}^{\infty} a^L q^{L^2} \alpha_L=
\sum_{L=0}^{\infty} a^L q^{L^2} \beta_L.
\end{equation}
The proof of the Rogers--Ramanujan and many similar such $q$-series identities
requires the input of suitable Bailey pairs into \eqref{ab}.
For example, from Rogers' work \cite{Rogers17} one can infer the 
following Bailey pair relative to $1$: $\alpha_0=1$ and
\begin{align}\label{BPRogers}
\alpha_L=(-1)^L q^{L(3L-1)/2}(1+q^L),\qquad
\beta_L=\frac{1}{(q)_L}.
\end{align}
Thus one finds
\begin{equation}\label{RR}
\frac{1}{(q)_{\infty}}\sum_{L=-\infty}^{\infty}
(-1)^L q^{L(5L-1)/2}=\sum_{n=0}^{\infty} \frac{q^{n^2}}{(q)_n}.
\end{equation}
The application of the Jacobi triple product identity
\begin{equation}\label{tpi}
\sum_{r=-\infty}^{\infty}(-1)^r x^r q^{\binom{r}{2}}=(x,q/x,q)_{\infty}
\end{equation}
yields the first Rogers--Ramanujan identity
\begin{equation*}
\sum_{n=0}^{\infty} \frac{q^{n^2}}{(q)_n}=
\frac{1}{(q,q^4;q^5)_{\infty}}.
\end{equation*}
Here and later in the paper we employ the condensed notation
$$(a_1,\dots,a_k;q)_n=(a_1,\dots,a_k)_n=(a_1)_n\dots (a_k)_n.$$
The second Rogers--Ramanujan identity 
\begin{equation*}
\sum_{n=0}^{\infty} \frac{q^{n(n+1)}}{(q)_n}=
\frac{1}{(q^2,q^3;q^5)_{\infty}}
\end{equation*}
follows in a similar fashion using the Bailey pair \cite{Rogers17} 
\begin{align*}
\alpha_L=(-1)^L q^{L(3L+1)/2}(1-q^{2L+1})/(1-q),\qquad
\beta_L=\frac{1}{(q)_L}
\end{align*}
relative to $q$.
By collecting as many Bailey pairs as possible, Slater compiled a list of
over a hundred Rogers--Ramanujan-type identities \cite{Slater51,Slater52}.
(Apart from a few exceptions Slater either used \eqref{ab} or the analogous 
identity obtained from \eqref{gd} and \eqref{abcd} by taking 
$M,\rho_1\to\infty$ and letting $\rho_2=-q^{k/2}$ with $k$ a small 
nonnegative integer.)

By dismissing the conjugate Bailey pair \eqref{gd} Bailey missed a very 
powerful mechanism for generating Bailey pairs. Namely, if we substitute 
the conjugate pair \eqref{gd} into \eqref{abcd} the resulting equation has 
the same form as the defining relation \eqref{BP} of a 
Bailey pair. This is formalized in the following theorem due to 
Andrews~\cite{Andrews84,Andrews85}.
\begin{theorem}\label{AB}
Let $(\alpha,\beta)$ form a Bailey pair relative to $a$. Then
so does $(\alpha',\beta')$ with
\begin{equation}\label{iter}
\begin{split}
\alpha_L'&=\frac{(\rho_1)_L(\rho_2)_L(aq/\rho_1\rho_2)^L}
{(aq/\rho_1)_L(aq/\rho_2)_L} \alpha_L \\[2mm]
\beta_L'&=\sum_{r=0}^L \frac{(\rho_1)_r(\rho_2)_r(aq/\rho_1\rho_2)^r
(aq/\rho_1\rho_2)_{L-r}}{(aq/\rho_1)_L(aq/\rho_2)_L(q)_{L-r}}\beta_r.
\end{split}
\end{equation}
\end{theorem}
Again letting $\rho_1,\rho_2$ tend to infinity leads to the important special
case
\begin{equation}\label{iter2}
\alpha_L'=a^L q^{L^2}\alpha_L \quad \text{and}\quad
\beta_L'=\sum_{r=0}^L \frac{a^r q^{r^2}}{(q)_{L-r}}\beta_r,
\end{equation}
which was also discovered by Paule~\cite{Paule85} for $a=1$ and $a=q$.

With this last result one finds that the Bailey pair of equation
\eqref{BPRogers} can be obtained by application of Theorem~\ref{AB}
with initial Bailey pair $\alpha_0=1$ and
\begin{equation}\label{initial}
\alpha_L=(-1)^L q^{\binom{L}{2}}(1+q^L),\qquad
\beta_L=\delta_{L,0}
\end{equation}
relative to $1$. Here $\delta_{i,j}$ is the Kronecker-delta symbol.
The Bailey pair \eqref{initial} follows after setting $x=1$ in the
$q$-binomial sum 
\begin{equation*}
\sum_{r=-L}^L (-1)^r x^r q^{\binom{r}{2}}\qbin{2L}{L-r}=(x,q/x)_L,
\end{equation*}
where throughout this paper the following definition of the
$q$-binomial coefficient or Gaussian polynomial is used
\begin{equation*}
\qbin{n}{m}=\frac{(q^{n-m+1})_m}{(q)_m}
\end{equation*}
for $m\in\Zp$ and zero otherwise.

At this point one may wonder why Bailey and Slater put so much emphasis
on finding new Bailey pairs, but contented themselves with just the single
conjugate Bailey pair~\eqref{CBP}. After all, the defining relations 
\eqref{BP} and \eqref{CBP} are very similar and it is therefore not 
unreasonable to expect that conjugate pairs are as important and as
numerous as ordinary Bailey pairs. In fact, because Andrews' Theorem~\ref{AB}
is equivalent to equation \eqref{abcd} with conjugate Bailey pair 
\eqref{gd}, in modern expositions of the Bailey lemma there often is
no mention of conjugate Bailey pairs and equation \eqref{abcd}
at all, see e.g. Refs.~\cite{Andrews84,Andrews86a,Andrews86b,Andrews92,
AAB87,ADH88,Bressoud88,BMS96,FQ96,Paule87,Paule88}.
Instead, \eqref{iter} is referred to as the Bailey lemma and in the
spirit of Slater, all focus is on finding interesting Bailey pairs.
These are then either iterated using \eqref{iter} or \eqref{iter2} to yield
what is called a Bailey chain, or directly substituted into \eqref{ab}.
The only exception that we were able to trace in the
literature is the conjugate Bailey pair ($|r|<1$)
\begin{equation}\label{BS}
\gamma_L=\frac{r^L}{(r)_{\infty}(aq)_{\infty}}
\sum_{k=0}^{\infty}\frac{(-1)^k q^{\binom{k}{2}}(r)_k(aq^{2L+1})^k}{(q)_k}
\quad \text{and}\quad \delta_L=r^L
\end{equation}
which can be found in the work by Bressoud~\cite{Bressoud81}
and Singh~\cite{Singh94} (and which for $r=q$ we will meet again in
Section~\ref{secSC}).

This paper intends to revive the interest in conjugate Bailey pairs.
In our earlier papers~\cite{SW97,SW98} we made a first step towards
this goal by proving an infinite series of conjugate Bailey pairs
generalizing \eqref{gdinf}.
Here we develop the theory of conjugate Bailey pairs much further,
exploiting the connection of Bailey's lemma with 
integrable systems and Lie algebras.
We show that appropriate series of one-dimensional configuration sums 
and A$_1^{(1)}$ string functions can be
identified with the series $\delta$ and $\gamma$ defining a conjugate
Bailey pair.
Here one-dimensional configuration sums~\cite{ABF84,FB85}, also known as
hook-partition generating functions~\cite{ABBBFV87}, are polynomials
that have arisen in statistical mechanics and partition theory.
A well-known example are the polynomials
introduced by Schur \cite{Schur17} in his famous proof of 
the Rogers--Ramanujan identities.
The string functions that occur are associated to the admissible
representations of the affine Kac--Moody algebra A$_1^{(1)}$ 
as introduced by Kac and Wakimoto~\cite{KW88,KW89}. 

Before we carry out the above program let us attempt to give an explanation
of the origin of our findings.
An important notion in the theory of affine Lie algebras is that of
branching functions~\cite{Kac90}.
Here we consider the branching functions $B^{N_1,N_2}$ associated to 
$(\text{A}_1^{(1)}\oplus \text{A}_1^{(1)},\text{A}_1^{(1)})$ at levels
$N_1$, $N_2$ and $N_1+N_2$, respectively, where $N_1$ and $N_2$
are rational numbers such that either $N_1$ or $N_2$ is a positive integer.
The branching functions obey the symmetry $B^{N_1,N_2}=B^{N_2,N_1}$.
Following the work of Andrews~\cite{Andrews84} and Foda and Quano~\cite{FQ96},
the infinite hierarchy of conjugate Bailey pairs of \cite{SW97,SW98}
were used in \cite{BMSW97} to derive $q$-series identities for the 
A$_1^{(1)}$ branching functions. Schematically the results of
\cite{BMSW97} read as follows:
\begin{equation}\label{BF12}
B^{N_1,N_2}=\sum_{L=0}^{\infty} \alpha_L^{(N_1)}\gamma_L^{(N_2)}
=\sum_{L=0}^{\infty} \beta_L^{(N_1)}\delta_L^{(N_2)},
\end{equation}
where $N_1$ is rational, $N_2$ integer, $\gamma_L^{(N_2)}$ is
a level-$N_2$ string function, $\beta_L^{(N_1)}$ a (normalized)
one-dimensional configuration sum and $(\alpha^{(N_1)},\beta^{(N_1)})$, 
$(\gamma^{(N_2)},\delta^{(N_2)})$ are a Bailey and conjugate Bailey
pair respectively.
In the middle term of this identity the symmetry between $N_1$ and
$N_2$ is not at all manifest since it involves only the (integer)
level-$N_2$ string functions and not the (fractional) level-$N_1$
string functions.  
This suggests that there should be more general conjugate Bailey
pairs such that one can also derive 
\begin{equation}\label{BF21}
B^{N_1,N_2}=\sum_{L=0}^{\infty} 
\bar{\alpha}_L^{(N_2)}\bar{\gamma}_L^{(N_1)}
=\sum_{L=0}^{\infty} \bar{\beta}_L^{(N_2)}\bar{\delta}_L^{(N_1)},
\end{equation}
where now $\bar{\gamma}_L^{(N_1)}$ is a fractional-level string function,
$(\bar{\alpha}^{(N_2)},\bar{\beta}^{(N_2)})$, 
$(\bar{\gamma}^{(N_1)},\bar{\delta}^{(N_1)})$ are a Bailey and conjugate 
Bailey pair, and such that, manifestly,
$$\sum_{L=0}^{\infty} \beta_L^{(N_1)}\delta_L^{(N_2)}
=\sum_{L=0}^{\infty} \bar{\beta}_L^{(N_2)}\bar{\delta}_L^{(N_1)}. $$
This last equation is obviously satisfied if
\begin{equation}\label{duality}
\bar{\delta}_L^{(N_1)}=g_L \beta_L^{(N_1)}\qquad\text{and}\qquad
\delta_L^{(N_2)}=g_L \bar{\beta}_L^{(N_2)},
\end{equation}
with $g_L$ independent of $N_1$ and $N_2$.
Since $\beta_L^{(N_1)}$ is a (normalized) one-dimensional configuration sum
we can now conclude that in the ``yet to be found'' conjugate Bailey pair
$(\bar{\gamma}^{(N_1)},\bar{\delta}^{(N_1)})$ the
sequence $\bar{\gamma}^{(N_1)}$ is a sequence
of (fractional) level-$N_1$ string functions and the sequence
$\bar{\delta}^{(N_1)}$ is proportional to a sequence
of one-dimensional configuration sums. This is indeed in accordance with
the announced results. We also note that the above discussion
establishes a duality between Bailey and conjugate Bailey pairs
through equation \eqref{duality}.

The remainder of the paper can be outlined as follows.
In the next two sections we review the one-dimensional configuration 
sums of the Andrews--Baxter--Forrester models and the
string functions associated with admissible representations
of the affine Lie algebra A$_1^{(1)}$. In Section~\ref{sec fb} these
are used to prove a very general class of conjugate Bailey pairs stated
in Corollary~\ref{cor CBP}.
In Section~\ref{sec fermi} we give fermionic or constant-sign expressions
for the one-dimensional configuration sums.
This allows us to apply the Bailey lemma,
together with our new conjugate Bailey pairs, to derive many new
$q$-series results in Sections~\ref{secSC} and \ref{sec BF}. 
In Section~\ref{secSC} we give fermionic formulas for the fractional-level 
A$_1^{(1)}$ string functions and parafermion characters.
In Section~\ref{sec BF} we derive a new type of bose-fermi identities
extending identities of the form \eqref{BF12} and \eqref{BF21}
for the A$_1^{(1)}$ branching functions by allowing for both
$N_1$ and $N_2$ to be rational numbers. To put this in the right context we
first present a discussion of the A$_1^{(1)}$ branching functions
in Section~\ref{sec branchingf} proving a generalization of a theorem of
Kac and Wakimoto that expresses the branching functions in terms
of fractional-level string functions in accordance with \eqref{BF21}.

\section{One-dimensional configuration sums}\label{sec cs}
The one-dimensional configuration sums of the
Andrews--Baxter--Forrester models were introduced in several
stages in Refs.~\cite{Schur17,Andrews81,ABF84,FB85}.
\begin{definition}\label{def cs}
For integers $p,p'$ with $1\leq p<p'$, and $b,s\in\N_{p'-1}$, $r\in\Z_{p+1}$
and $L\in\Zp$ such that $L+s+b$ is even, let
\begin{multline}\label{odcs}
X^{(p,p')}_{r,s}(L,b;q)=
X^{(p,p')}_{r,s}(L,b) \\[2mm]
=\sum_{j\in\Z} \Bigr\{q^{j(p p'j+p'r-ps)}\qbins{L}{(L+s-b)/2-p'j}-
q^{(pj+r)(p'j+s)}\qbins{L}{(L-s-b)/2-p'j}\Bigl\}.
\end{multline}
\end{definition}
The configuration sums possess two symmetries which will be used later.
{}From the definition it can be deduced immediately that 
\begin{equation}\label{flip symmetry}
X_{r,s}^{(p,p')}(L,b) = X_{p-r,p'-s}^{(p,p')}(L,p'-b),
\end{equation}
whereas 
\begin{equation}\label{dual symmetry}
X_{r,s}^{(p,p')}(L,b;q)=
q^{\frac{1}{4}(L^2-(b-s)^2)} X_{b-r,s}^{(p'-p,p')}(L,b;1/q)
\end{equation}
follows by application of 
\begin{equation}\label{qbininv}
\qbin{n}{m}_{1/q}=q^{m(m-n)}\qbin{n}{m}.
\end{equation}

When the parameters $p$ and $p'$ obey the additional restriction
\begin{equation}\label{cond}
\gcd(p,p')=1
\end{equation}
the polynomials \eqref{odcs} were encountered by Forrester and 
Baxter~\cite{FB85} as the generating function of sets of restricted 
lattice path. Below we describe a slight extension of their result.
A lattice path interpretation of the one-dimensional configuration sums 
$X^{(p,p')}_{r,s}(L,b;q)$ for all $1\leq p<p'$ can be found in~\cite{FLPW98}.

Let $P=(x_0,\dots,x_{L+1})$ be a lattice path consisting of an ordered
sequence of $L+2$ integers such that
$|x_{i+1}-x_i|=1$ for $0\leq i\leq L$, $x_0=s$, $x_L=b$, $x_{L+1}=c$ and
$x_i\in\N_{p'-1}$ for $1\leq i\leq L$. Denote the set of all such paths by 
$\P_L^{s,b,c}$.
Assign a weight $|P|$ to $P\in\P_L^{s,b,c}$ as follows
\begin{equation*}
|P|=\sum_{i=1}^L i H(x_{i-1},x_i,x_{i+1}),
\end{equation*}
where
\begin{equation*}
H(a,a\mp 1,a)=\pm \Big\lfloor \frac{a(p'-p)}{p'}\Big\rfloor 
\quad\text{and}\quad H(a\pm 1,a,a\mp 1)=\frac{1}{2}.
\end{equation*}
Here $\lfloor x \rfloor$ denotes the integer part of $x$.
Forrester and Baxter studied the generating function 
\begin{equation*}
D_L(s,b,c;q)=\sum_{P\in\P_L^{s,b,c}}q^{|P|}
\end{equation*}
and proved for $c\in\N_{p'-1}$ that \cite[Thm 2.3.1]{FB85}
\begin{equation}\label{DX}
D_L(s,b,c;q)=q^{\frac{1}{4}L(c-b)(c+b-1-2r)+\frac{1}{4}(s-b)(s+b-1-2r)}
X_{r,s}^{(p,p')}(L,b),
\end{equation}
where $r$ is given by
\begin{align}\label{kr}
r&=\frac{b+c-1}{2}-\Big\lfloor\frac{c(p'-p)}{p'}\Big\rfloor \\
\label{kr2} &=\frac{b-c+1}{2}+\Big\lfloor\frac{cp}{p'}\Big\rfloor .
\end{align}
For $p'=p+1$ this result was first obtained in \cite{ABF84}.

Later in this paper the configuration sum $X_{0,s}^{(p,p')}(L,1)$
will play a prominent role. Using the standard $q$-binomial recurrences
\begin{equation*}
\qbin{n}{m}=\qbin{n-1}{m-1}+q^m \qbin{n-1}{m}=
\qbin{n-1}{m}+q^{n-m}\qbin{n-1}{m-1}
\end{equation*}
it readily follows that
\begin{equation*}
\qbin{L}{a}-\qbin{L}{a-1}= q^a\qbin{L}{a}-q^{L-a+1}\qbin{L}{a-1}.
\end{equation*}
One thus finds the relation
\begin{equation}\label{rzero}
X_{0,s}^{(p,p')}(L,1) = q^{\frac{1}{2}(L-s+1)} X_{1,s}^{(p,p')}(L,1).
\end{equation}
The corresponding lattice path interpretation for $X_{0,s}^{(p,p')}(L,1)$
is easily found. When $p'>2p$ it is included in the Forrester--Baxter result
\eqref{DX} since $b=1$ and $c=2$ yields $r=0$.
When $p'<2p$ we need to allow for paths with $c=0$. 
Then $b=1$  and, using \eqref{kr}, $r=0$.
To see that the corresponding generating function is indeed
\begin{equation}\label{D}
D_L(s,1,0;q)=q^{\frac{1}{4}s(s-1)}X_{0,s}^{(p,p')}(L,1)
\end{equation}
we compute $D_L(s,1,2;q)/D_L(s,1,0;q)$.
On the one hand, by the one-to-one correspondence 
$(s,x_2,\dots,x_{L-2},2,1,2)\leftrightarrow (s,x_2,\dots,x_{L-2},2,1,0)$
between paths in $\P_L^{s,1,2}$ and $\P_L^{s,1,0}$,
and the fact that  $H(2,1,2)=0$ and $H(2,1,0)=1/2$ one finds
$D_L(s,1,2;q)/D_L(s,1,0;q)=q^{-L/2}$.
On the other hand, by \eqref{DX} and \eqref{rzero} we get
\begin{equation*}
\frac{D_L(s,1,2;q)}{D_L(s,1,0;q)}=q^{\frac{1}{4}(s-1)(s-2)} 
\frac{X_{1,s}^{(p,p')}(L,1)}{D_L(s,1,0;q)}
=q^{\frac{1}{4}s(s-1)-\frac{1}{2}L} 
\frac{X_{0,s}^{(p,p')}(L,1)}{D_L(s,1,0;q)}.
\end{equation*}
Combining the last two results clearly implies \eqref{D}.

By the symmetry \eqref{flip symmetry} we also need $X_{p,s}^{(p,p')}(L,p'-1)$.
For $p'>2p$ its lattice path interpretation follows again from the
Forrester--Baxter result, as $b=p'-1$ and $c=p'-2$ yields $r=p$.
When $p'<2p$ we need to allow for paths with $c=p'$.
Then $b=p'-1$ and, using \eqref{kr2}, $r=p$.
By a calculation similar to the one above it is then readily shown that
$D_L(s,p'-1,p';q)$ is indeed given by \eqref{DX}.

The expressions \eqref{odcs} have also been studied extensively in the
theory of partitions, see e.g. 
\cite{Andrews72,ABBBFV87,Bressoud80a,Bressoud96,Burge82,GK97}.
Here we quote the most general result, obtained in \cite{ABBBFV87}.
Let $\lambda$ be a partition and $\lambda'$ its conjugate.
The $(i,j)$-th node of $\lambda$ is the node (or box)
in the $i$th row and $j$th column of the Ferrers diagram of $\lambda$.
The $d$th diagonal of $\lambda$ is formed by the nodes 
with coordinates $(i,i-d)$. The hook difference at node $(i,j)$ is defined
as $\lambda_i-\lambda_j'$. Theorem~1 of \cite{ABBBFV87} states that
the generating function of partitions $\lambda$
with at most $(L+s-b)/2$ parts, largest part not exceeding $(L-s+b)/2$, 
and hook differences on the $(1-r)$th diagonal at least $r-s+1$ and on 
the $(p-r-1)$th diagonal at most $p'-p+r-s-1$ is given by 
$X_{r,s}^{(p,p')}(L,b)$. Here the following two conditions 
apply~\cite{ABBBFV87}, $1\leq r\leq p-1$ and $0\leq b-r\leq p'-p$.
When $r=0$ one has to impose the additional condition that the largest
part exceeds $(L-s-b)/2$.
Similarly, the case $r=p$ can be included provided one demands that the
number of parts exceeds $(L+s+b)/2$.

\section{Characters and string functions for A$_1^{(1)}$}\label{sec string}

In~\cite{KW88,KW89} Kac and Wakimoto introduced admissible highest weight 
representations of affine Lie algebras as generalizations of
the familiar integrable highest weight representations~\cite{Kac90}.
Let $p,p'$ be integers such that $1\leq p<p'$ and $\gcd(p,p')=1$, and define
$$N=p'/p-2$$ so that $-1<N<0$ for $p<p'<2p$ and $N>0$ for $p'>2p$.
Let $\Lambda_0$ and $\Lambda_1$ be the fundamental weights of A$_1^{(1)}$.
Fix an integer $\ell\in\Z_{p'-1}$ and let $L(\lambda)$ be an admissible 
A$_1^{(1)}$ highest weight module of highest weight\footnote{Kac and 
Wakimoto considered the more general case
$\lambda=(N-\ell)\Lambda_0+\ell\Lambda_1+k(N+2)(\Lambda_0-\Lambda_1)$
with $k\in\Z_p$.}
$\lambda=(N-\ell)\Lambda_0+\ell\Lambda_1$.
The corresponding character is formally defined as
\begin{equation*}
\chi_{\ell}^N(z,q)=\chi_{\ell}(z,q)=\text{tr}_{L(\lambda)}
q^{s_{\lambda}-d} z^{-\frac{1}{2}\alpha_1^\vee},
\end{equation*}
where $d=3$ is the dimension of A$_1$, $\alpha_1^\vee$ is a simple coroot
and $$s_{\lambda}=-\frac{1}{8}+\frac{(\ell+1)^2}{4(N+2)}.$$
In terms of the classical theta function 
\begin{equation}\label{Thetadef}
\Theta_{n,m}(z,q)=\sum_{j\in\Z+n/2m} q^{mj^2}z^{-mj}
\end{equation}
of degree $m$ and characteristic $n$, one can express the A$_1^{(1)}$
character as
\begin{equation}\label{char}
\chi_{\ell}(z,q)=\frac{\sum_{\sigma=\pm 1}\sigma
\Theta_{\sigma(\ell+1),p'}(z,q^p)}
{\sum_{\sigma=\pm 1}\sigma\Theta_{\sigma,2}(z,q)}.
\end{equation}

The level-$N$ A$_1^{(1)}$ string functions are defined by the expansion
\begin{equation}\label{sfdef}
\chi_{\ell}(z,q)=\sum_{m\in 2\Z+\ell} C_{m,\ell}^N(q) 
q^{\frac{m^2}{4N}} z^{-\frac{1}{2}m},
\end{equation}
and enjoy the symmetry
\begin{equation}\label{sym1}
C_{m,\ell}^N=C_{-m,\ell}^N.
\end{equation}
When $N$ is integer we furthermore have
\begin{equation}\label{sym2}
C_{m,\ell}^N=C_{2N-m,\ell}^N=C_{N-m,N-\ell}^N
\end{equation}
so that \eqref{sfdef} may be put in the familiar form
\begin{equation*}
\chi_{\ell}(z,q)=\sum_{\substack{0\leq m <2N\\ \text{$m+\ell$ even}}}
C_{m,\ell}^N(q)\Theta_{m,N}(z,q).
\end{equation*}

We derive an expression for the string functions following the
approach of e.g., Refs.~\cite{HNY90,ACT91}. First observe that
\begin{equation*}
\sum_{\sigma=\pm 1}\sigma\Theta_{\sigma,2}(z,q)=
q^{1/8}z^{-\frac{1}{2}}\sum_{j\in\Z}(-1)^j q^{\binom{j}{2}}z^j=
q^{1/8}z^{-\frac{1}{2}}(z,q/z,q)_{\infty}
\end{equation*}
where in the second step Jacobi's triple product identity \eqref{tpi} 
has been employed. Next recall the identity 
\begin{equation*}
\frac{1}{(z,q/z)_{\infty}}=
\frac{1}{(q)_{\infty}^2}
\sum_{k\in\Z}\sum_{i\in\N}
(-1)^{i+1}q^{\binom{i}{2}-ik}z^k,
\end{equation*}
which can be extracted from an expansion of the following ratio of 
Jacobi theta functions $\vartheta_1'(0)/\vartheta_1(u)$
in {\cite[\S 486]{TM98} (see also \cite[Eq. (5.26)]{KP84} and 
\cite[Eqs. (A.4), (A.5)]{Thorn89}). Using this we find that
\begin{multline*}
\chi_{\ell}(z,q)=\frac{1}{\eta^3(\tau)}
\sum_{\sigma=\pm 1}\sum_{j,k\in\Z}\sum_{i\in\N}
\sigma (-1)^{i+1}q^{\binom{i}{2}-ik+
pp'(j+\sigma(\ell+1)/(2p'))^2} \\
\times z^{-\frac{1}{2}(2p'j-2k+\sigma(\ell+1)-1)},
\end{multline*}
where, as usual, $\eta(\tau)=q^{1/24}(q)_{\infty}$ with
$q=\exp(2\pi i\tau)$.
Now replace $j$ by $\sigma j$ and then 
$k$ by $\frac{1}{2}(2\sigma p'j-m-1+\sigma(\ell+1))$. This yields
\begin{multline*}
\chi_{\ell}(z,q)=\frac{1}{\eta^3(\tau)}
\sum_{m\in 2\Z+\ell} \sum_{j\in\Z}\sum_{i\in\N}
(-1)^i q^{\frac{1}{2}i(i+m)+ pp'(j+(\ell+1)/(2p'))^2} \\
\times \Bigl\{q^{\frac{1}{2}i(2p'j+\ell+1)}
-q^{-\frac{1}{2}i(2p'j+\ell+1)}\Bigr\}z^{-\frac{1}{2}m}.
\end{multline*}
Comparing this with \eqref{sfdef} one can extract the string functions as
\begin{multline}\label{sffinal}
C_{m,\ell}^N(q)=\frac{q^{\frac{(\ell+1)^2}{4(N+2)}-\frac{m^2}{4N}}}
{\eta^3(\tau)}\sum_{j\in\Z}
\sum_{i\in\N} (-1)^i q^{\frac{1}{2}i(i+m)+jp(p'j+\ell+1)} \\
\times
\Bigl\{q^{\frac{1}{2}i(2p'j+\ell+1)}-q^{-\frac{1}{2}i(2p'j+\ell+1)}\Bigr\}.
\end{multline}

We slightly extend the original definition of the string functions
given in equation~\eqref{sfdef} by dropping the condition $\gcd(p,p')=1$.
Also normalizing for later convenience we are led to the following
definition.
\begin{definition}
For integers $1\leq p<p'$, $m\in\Z$ and $\ell\in\Z_{p'-1}$
such that $\ell$ and $m$ have equal parity,
\begin{multline}\label{sfdef2}
\C_{m,\ell}^{(p,p')}(q)=\\
\frac{1}{(q)_{\infty}^3}\sum_{j\in\Z}
\sum_{i\in\N} (-1)^i q^{\frac{1}{2}i(i+m)+jp(p'j+\ell+1)}
\Bigl\{q^{\frac{1}{2}i(2p'j+\ell+1)}-q^{-\frac{1}{2}i(2p'j+\ell+1)}\Bigr\}.
\end{multline}
\end{definition}
When $\gcd(p,p')=1$ we also use the notation 
$\C_{m,\ell}^{N}(q)=\C_{m,\ell}^{(p,p')}(q)$,
where $N=p'/p-2$ is the level of the modified string function.

As a note of warning we remark that for a generic choice of variables
the order of summation in \eqref{sffinal} and \eqref{sfdef2} has to
be strictly obeyed. We use the form \eqref{sfdef2} as defining relation
rather than the more familiar (and computationally efficient) expression
\begin{align}\label{sfstand}
\C_{m,\ell}^{(p,p')}(q)&=\frac{1}{(q)_{\infty}^3}
\Bigl\{\sum_{\substack{i\geq 0\\[.5mm] j\geq 0}}-
\sum_{\substack{i<0 \\[.5mm]j<0}}\Bigr\}
(-1)^i q^{\frac{1}{2}i(i+m)+p'j(pj+i)+\frac{1}{2}(\ell+1)(2pj+i)} \\
&-\frac{1}{(q)_{\infty}^3}
\Bigl\{\sum_{\substack{i\geq 0 \\[.5mm]j>0}}-
\sum_{\substack{i<0\\[.5mm] j\leq 0}}\Bigr\}
(-1)^i q^{\frac{1}{2}i(i+m)+p'j(pj+i)-\frac{1}{2}(\ell+1)(2pj+i)} \notag
\end{align}
for later reasons. By
\begin{equation}\label{zero}
\sum_{i=-\infty}^{\infty} (-1)^i q^{\binom{i}{2}+in}=0\quad 
\text{for $n\in\Z$},
\end{equation}
which is a specialization of Jacobi's triple product identity \eqref{tpi},
it is straightforward to transform \eqref{sfdef2} into \eqref{sfstand}.
We also note that for integer level, i.e., $p=1$ and $p'=N+2$ we can rewrite
\eqref{sfstand} in the neat form (by \eqref{zero} equivalent 
to~\cite[Eq. (3.17)]{DQ90}) 
\begin{equation*}
\C_{m,\ell}^N(q)=\frac{1}{(q)_{\infty}^3}
\Bigl\{\sum_{\substack{j\geq 1\\[.5mm] k\leq 0}}
-\sum_{\substack{j\leq 0\\[.5mm] k\geq 1}}\Bigr\}
(-1)^{k-j}q^{\binom{k-j}{2}-Njk+\frac{1}{2}k(m-\ell)+\frac{1}{2}j(m+\ell)}.
\end{equation*}
To see this, make the variable changes $j\to -j$ followed by $i\to k+j-1$
in the first line and $j\to 1-k$ followed by $i\to k+j-1$ in the second
line of \eqref{sfstand} and use the symmetry
$\C_{m,\ell}^N(q)=q^{(m-\ell)/2}\C_{m-N,N-\ell}^N(q)$.

To conclude this section we introduce the characters 
$e_{m,\ell}^N(q)$ of the Z$_N$ para\-fer\-mi\-on algebra at rational 
level $N$~\cite{ACT91}. It was argued in \cite{ACT91} that these
characters are realized as branching functions as follows:
\begin{equation*}
\chi_{\ell}(z,q)=\sum_{m\in 2\Z+\ell} 
e_{m,\ell}^N(q) \frac{q^{\frac{m^2}{4N}}z^{-m/2}}{\eta(\tau)}
\end{equation*}
Comparison with \eqref{sfdef} shows that
\begin{equation}\label{edef}
e_{\ell,m}^N(q)=\eta(\tau) C_{\ell,m}^N(q).
\end{equation}
For integer $N$ the $e_{m,\ell}^N$ have also been shown to be
branching functions of the Lie algebra pair (A$^{(1)}_{2N-1}$, C$^{(1)}_N$)
\cite{JM84}.

\section{Fractional-level conjugate Bailey pairs}\label{sec fb}
This section contains the key results of this paper.
In Theorem~\ref{thmCBP} new conjugate Bailey pairs are stated, 
which by Corollary~\ref{cor CBP} imply conjugate Bailey pairs involving
the one-dimensional configurations sums and fractional-level string functions
of the previous two sections.

\begin{theorem}\label{thmCBP}
For $\eta\in\Zp$ and $j\in\Z$, the pair of sequences $(\gamma,\delta)$
with 
\begin{equation}\label{CBP1}
\begin{split}
\gamma_L&=\frac{1}{(q)_{\infty}^2(aq)_{\infty}}
\sum_{i=1}^{\infty}(-1)^i q^{\frac{1}{2}i(i+2L+\eta)}
\Bigl\{q^{\frac{1}{2}i(2j+\eta+1)}-q^{-\frac{1}{2}i(2j+\eta+1)}\Bigr\}\\[1mm]
\delta_L&=\qbins{2L+\eta}{L-j}-\qbins{2L+\eta}{L-j-1}
\end{split}
\end{equation}
forms a conjugate Bailey pair relative to $a=q^{\eta}$.
\end{theorem}

Before we prove this theorem let us first state the following corollary.
\begin{corollary}\label{cor CBP}
Fix integers $1\leq p<p'$, and
let $\eta\in\Zp$ and $\ell\in\Z_{p'-1}$ such that
$\ell+\eta$ is even. Let $\C_{m,\ell}^{(p,p')}$ and
$X^{(p,p')}_{r,s}(L,b)$ be defined as in \eqref{sfdef2} and \eqref{odcs}.
Then $(\gamma,\delta)$ with 
\begin{equation}\label{CBP2}
\gamma_L=(q)_{\eta}\C_{2L+\eta,\ell}^{(p,p')}(q)  \qquad \text{and} \qquad
\delta_L=X^{(p,p')}_{0,\ell+1}(2L+\eta,1)
\end{equation}
forms a conjugate Bailey pair relative to $a=q^{\eta}$.
\end{corollary}

\begin{proof}
Take the conjugate Bailey pair \eqref{CBP1} and replace
$j$ by $jp'+(\ell-\eta)/2$. Then multiply both $\gamma_L$ and
$\delta_L$ by $q^{jp(jp'+\ell+1)}$ and sum $j$ over the integers.
Using \eqref{sfdef2} and \eqref{odcs} this transforms $\gamma_L$ 
and $\delta_L$ of \eqref{CBP1} into those of \eqref{CBP2}.
\end{proof}

The proof of Theorem~\ref{thmCBP} rests upon the following lemma.
\begin{lemma}\label{lemCBP}
For $a$ and $b$ indeterminates,
\begin{multline}\label{ablem}
\sum_{r=0}^{\infty}\frac{(ab)_{2r}}{(q)_r(ab)_r}\Bigl\{
\frac{1}{(aq)_{r-1}(bq)_r}-
\frac{1}{(aq)_r(bq)_{r-1}}\Bigr\}\\
=\frac{1}{(q)_{\infty}(aq)_{\infty}(bq)_{\infty}}
\sum_{i=1}^{\infty}(-1)^i q^{\binom{i}{2}}(a^i-b^i).
\end{multline}
\end{lemma}

\begin{proof}
The terms within the curly braces can be combined to
$(b-a)q^r/(aq)_r(bq)_r$. Using this as well as
$(a)_{\infty}/(a)_r=(aq^r)_{\infty}$ and $(a)_{2r}/(a)_r=(aq^r)_r$, 
equation~\eqref{ablem} can be written as
\begin{equation}\label{nf}
\sum_{r=0}^{\infty}
\frac{q^r (abq^r)_r(aq^{r+1})_{\infty}(bq^{r+1})_{\infty}}{(q)_r}
=\frac{1}{(q)_{\infty}}
\sum_{i=1}^{\infty}(-1)^{i+1} q^{\binom{i}{2}}\frac{a^i-b^i}{a-b}.
\end{equation}
We now use the $q$-binomial sum \cite[Eq. (3.3.6)]{Andrews76a}
\begin{equation}\label{qbinthm}
(a)_n=\sum_{k=0}^n (-a)^k q^{\binom{k}{2}}\qbin{n}{k}
\end{equation}
as well as the limiting case
\begin{equation*}
(a)_{\infty}=\sum_{k=0}^{\infty}\frac{(-a)^k q^{\binom{k}{2}}}{(q)_k},
\end{equation*}
to express the left-hand side of \eqref{nf} as the following quadruple sum,
\begin{equation*}
\sum_{r=0}^{\infty}\sum_{i=0}^{\infty}\sum_{j=0}^{\infty}\sum_{k=0}^r
(-1)^{i+j+k} a^{i+k}b^{j+k} \;
\frac{q^{\binom{i+1}{2}+\binom{j+1}{2}+\binom{k}{2}+r(i+j+k+1)}} 
{(q)_i(q)_j(q)_k(q)_{r-k}}.
\end{equation*}
After shifting $i\to i-k$, $j\to j-k$ and $r\to r+k$ this becomes
\begin{equation*}
\sum_{i=0}^{\infty}\sum_{j=0}^{\infty}(-1)^{i+j}
q^{\binom{i+1}{2}+\binom{j+1}{2}}a^i b^j 
\sum_{k=0}^{\min\{i,j\}}
\frac{(-1)^kq^{\binom{k}{2}}} 
{(q)_{i-k}(q)_{j-k}(q)_k}
\sum_{r=0}^{\infty}
\frac{q^{r(i+j-k+1)}}{(q)_r}.
\end{equation*}
The sum over $r$ can readily be performed thanks to
\cite[Eq. (1.3.15)]{GR90}
\begin{equation}\label{xq}
\sum_{r=0}^{\infty} \frac{x^r}{(q)_r}=\frac{1}{(x)_{\infty}},
\end{equation}
leading to
\begin{equation*}
\frac{1}{(q)_{\infty}}
\sum_{i=0}^{\infty}\sum_{j=0}^{\infty}(-1)^{i+j}
q^{\binom{i+1}{2}+\binom{j+1}{2}}a^i b^j 
\sum_{k=0}^{\min\{i,j\}}
\frac{(-1)^k q^{\binom{k}{2}}(q)_{i+j-k}}
{(q)_{i-k}(q)_{j-k}(q)_k}.
\end{equation*}
The sum over $k$ yields $q^{ij}$ by the $q$-Chu--Vandermonde sum
\cite[Eq. (II.7)]{GR90}
\begin{equation}\label{qCV}
{_2\phi_1}\Bigl[\genfrac{}{}{0pt}{}{a,q^{-n}}{c};q,\frac{cq^n}{a}\Bigr]
=\frac{(c/a)_n}{(c)_n}
\end{equation}
with $n=\min\{i,j\}$, $c=q^{-i-j}$ and $a=cq^n$,
where the following standard notation for basic hypergeometric 
series is employed
\begin{equation*}
{_{r+1}\phi_r}\Bigl[\genfrac{}{}{0pt}{}{a_1,\dots,a_{r+1}}
{b_1,\dots,b_r};q,z\Bigr]
=\sum_{k=0}^{\infty}\frac{(a_1,\dots,a_{r+1})_k}{(q,b_1,\dots,b_r)_k}z^k.
\end{equation*}
As a result we are left with
\begin{equation*}
\frac{1}{(q)_{\infty}}
\sum_{i=0}^{\infty}\sum_{j=0}^{\infty}(-1)^{i+j}a^i b^j 
q^{\binom{i+j+1}{2}}.
\end{equation*}
This corresponds to the right-hand side of \eqref{nf} as
\begin{multline*}
\sum_{i=1}^{\infty}(-1)^{i+1} q^{\binom{i}{2}}\frac{a^i-b^i}{a-b}
=\sum_{i=1}^{\infty}(-1)^{i+1} q^{\binom{i}{2}}
\sum_{j=0}^{i-1} a^{i-j-1} b^j \\
=\sum_{j=0}^{\infty}\sum_{i=j+1}^{\infty} 
(-1)^{i+1} q^{\binom{i}{2}}a^{i-j-1} b^j 
=\sum_{j=0}^{\infty}\sum_{i=0}^{\infty} 
(-1)^{i+j} q^{\binom{i+j+1}{2}} a^i b^j.
\end{multline*}
\end{proof}
Finally we have to show that Theorem~\ref{thmCBP} follows from
Lemma~\ref{lemCBP}. 
\begin{proof}[Proof of Theorem~\ref{thmCBP}]
Substitute the conjugate Bailey pair \eqref{CBP1} into the defining
relation \eqref{CBP}. After the shift $r\to r+L$ this becomes
\begin{multline*}
\frac{1}{(q)_{\infty}^3} \sum_{i=1}^{\infty}(-1)^i q^{\binom{i}{2}}
\Bigl\{q^{\frac{1}{2}i(\zeta-\sigma)}
-q^{\frac{1}{2}i(\zeta+\sigma+2)}\Bigr\}\\
=\sum_{r=0}^{\infty} \frac{1}{(q)_r(q)_{r+\zeta}}
\Bigl\{\qbins{2r+\zeta}{r+\frac{1}{2}(\zeta-\sigma-2)}-
\qbins{2r+\zeta}{r+\frac{1}{2}(\zeta-\sigma)}\Bigr\},
\end{multline*}
where we have set $2L+\eta=\zeta\geq 0$ and $2j+\eta=\sigma$.
To obtain this identity we take \eqref{ablem} and choose
$a=q^{(\zeta-\sigma)/2}$, $b=q^{(\zeta+\sigma+2)/2}$
and perform a few trivial operations.
\end{proof}

\section{Fermionic expressions for the one-dimensional configuration sums}
\label{sec fermi}

{}From Definition~\ref{def cs} of the one-dimensional configuration sums
we see that the sequence $\delta$ in Corollary~\ref{cor CBP} is not
a sequence of manifestly positive polynomials (polynomials with positive
integer coefficients). In applications of the corollary
interesting $q$-series identities arise when there exist expressions 
that do have this property.
Such constant-sign or fermionic representations for the configuration sums
of the Andrews--Baxter--Forrester models
have recently attracted a lot of attention~\cite{Melzer94,Berkovich94,FQ95,
BM96,Warnaar96a,Warnaar96b,BMS98,DF98,FLPW98,FW98}.
In this section we present some of the cited results for 
$X_{r,s}^{(p,p')}(L,b)$ in the simplest case when $\gcd(p,p')=1$ and 
$s$ and $b$ are so-called Takahashi lengths
associated with the continued fraction expansion of $p/(p'-p)$.
More complicated cases where $s$ and $b$ are not necessarily
Takahashi lengths or where $(p,p')\neq 1$
can be found in \cite{BMS98} and \cite{FQ95}, respectively.

Given $p,p'$ such that $\gcd(p,p')=1$ and $p<p'<2p$ define
integers $n$ and $\nu_0,\ldots,\nu_n$ by the continued fraction expansion
\begin{equation*}
\frac{p}{p'-p}=[\nu_0,\nu_1,\ldots,\nu_n].
\end{equation*}
Introduce partial sums of the $\nu_j$ as
$t_m=\sum_{j=0}^{m-1}\nu_j$ for $1\le m\le n$ and
set $t_0=-1$ and $d=t_{n+1}=t_n+\nu_n-2$.
The $t_m$'s define a matrix $\I_B$ of size $d\times d$ with entries
\begin{equation*}
(\I_B)_{i,j} = \begin{cases}
\delta_{i,j+1}+\delta_{i,j-1} & \text{for $i\neq t_m$}\\
\delta_{i,j+1}+\delta_{i,j}-\delta_{i,j-1} & 
\text{for $i=t_m<d$}\\
\delta_{i,j+1}+\delta_{\nu_n,2}\delta_{i,j} & \text{for $i=d$.}
\end{cases}
\end{equation*}
Viewing $\I_B$ as a generalized incidence matrix we define a
corresponding fractional-level Cartan-type matrix $B=2I-\I_B$,
where $I$ is the identity matrix. 
When $p'=p+1$ the matrix $B$ is a Cartan matrix of type A
and when $p'=p+2$ it corresponds to a Cartan-type matrix of a tadpole graph.

For $1\le m\le n$ consider the recursion
\begin{equation*}
x_{m+1}=x_{m-1}+\nu_m x_m.
\end{equation*}
We need two sets of integers $\{y_m\}_{m=0}^{n+1}$ and 
$\{\bar{y}_m\}_{m=0}^{n+1}$ approximating $p'$ and $p$,
defined by the above recurrence and the initial conditions
$y_{-1}=0$, $\bar{y}_{-1}=-1$, $y_0=\bar{y}_0=1$
$y_1=\nu_0+1$, $\bar{y}_1=\nu_0$.
Hence $\bar{y}_m/(y_m-\bar{y}_m)=[\nu_0,\dots,\nu_{m-1}]$,
$y_{n+1}=p'$ and $\y_{n+1}=p$.
An important subset of $\N_{p'-1}$ is given
by the ``Takahashi lengths'' $l_1,\dots,l_{d+2}$ defined as
\begin{equation*}
l_{j+1} = y_{m-1}+(j-t_m)y_m, \qquad
t_m<j\le t_{m+1}+\delta_{m,n}.
\end{equation*}
Clearly, for $p'=p+1$ the set of Takahashi lengths is just $\N_{p'-1}$.
Similarly one may define the ``truncated Takahashi lengths''
$\lb_1,\dots,\lb_{d+2}$,
\begin{equation*}
\lb_{j+1} = \y_{m-1}+(j-t_m)\y_m, \qquad
t_m<j\le t_{m+1}+\delta_{m,n},
\end{equation*}
which determine a subset of $\Z_p$. If $b=l_{j+1}$ is a Takahashi length
then $\bar{b}$ denotes the truncated Takahashi length $\lb_{j+1}$.

For vectors $\vu,\vv\in\Zp^{d+1}$ define
\begin{equation}\label{f}
f(\vu,\vv)=\sum_{\vm\in 2\Z^d+\vQ_{\vu+\vv}}
q^{\frac{1}{4}\vm B\vm-\frac{1}{2}\vA_{\vu,\vv}\vm}\qbin{\vm+\vn}{\vm},
\end{equation}
where
\begin{equation*}
\qbin{\vm+\vn}{\vm}=\prod_{j=1}^d\qbin{m_j+n_j}{m_j}
\end{equation*}
and where the following definitions are used.
The variables $\vm$ and $\vn$ are related by the $(\vm,\vn)$-system
\begin{equation*}
\vm+\vn = \frac{1}{2}(\I_B\vm+\vu^*+\vv^*)
\end{equation*}
where $\vu^*$ and $\vv^*$ denote the projections of $\vu$ and $\vv$ 
onto $\Zp^d$. The linear term in the exponent of \eqref{f} is fixed by
\begin{equation*}
(\vA_{\vu,\vv})_k=\begin{cases}
u_k & \text{for $m$ odd}\\
v_k & \text{for $m$ even}
\end{cases} \qquad t_m<k\le t_{m+1}.
\end{equation*}
Finally, $\vQ_{\vu}=\sum_{j=1}^{d+1} u_j\vQ^{(j)}$
where $\vQ^{(j)}$ is defined recursively as
\begin{equation*}
Q_i^{(j)}=\begin{cases}
 \max\{j-i,0\} & \text{for $t_m\le i\le d$}\\[1mm]
Q_{i+1}^{(j)}+Q_{t_{m'}+1}^{(j)} & \text{for $t_{m'-1}\le i< t_{m'},
 1\le m'\le m$}
\end{cases}
\end{equation*}
with $0\le m \le n$ such that $t_m<j\le t_{m+1}
+\delta_{m,n}$. When $\nu_n=2$ we must take $Q_{t_n+1}^{(t_n+1)}=0$.

When the conditions \eqref{cond} are satisfied there exist
fermionic expressions for the one-dimensional 
configuration sums \eqref{odcs} in terms of the function
\eqref{f}~\cite{BMS98}.
Generally these are very complex and, as mentioned earlier, to keep formulas
relatively simple we restrict our attention to $b$ and $s$ being Takahashi
lengths (see \cite[Eq. (10.3)]{BMS98}).
\begin{theorem}
Let $1\le p<p'<2p$ such that $\gcd(p,p')=1$ and let
$b=l_{\beta+1}$, $s=l_{\sigma+1}$ be Takahashi lengths with $\beta\ge 1$
and $r=\bar{b}=\lb_{\beta+1}$. Then
\begin{equation}\label{F}
X_{r,s}^{(p,p')}(L,b) = 
q^{\Delta_{b,s}} f(L\ve_1+\vu_{\beta},\vu_{\sigma}),
\end{equation}
where $\ve_i$ is the $i$th standard unit vector in $\Z^{d+1}$ ($\ve_0=0$) and
\begin{equation}\label{vui}
\vu_i = \ve_i-\sum_{k=m+1}^n\ve_{t_k}\qquad
\text{for $t_m<i\le t_{m+1}+\delta_{m,n}$.}
\end{equation}
\end{theorem}
The explicit expression for $\Delta_{b,s}$ in the theorem is quite involved
and is omitted here. Instead we fix it by requiring that
\begin{equation*}
X_{r,s}^{(p,p')}(L,b;q=0)=1,
\end{equation*}
for $L\geq |s-b|$.
The relation between $b$ and $r$ given in the theorem 
corresponds to \eqref{kr} with $c=b-1$.
This explains why $\beta\geq 1$ (or $b=l_{\beta+1}\geq 2$).
As a consequence 
$X_{0,s}^{(p,p')}(L,1)$, or, equivalently,
$X_{1,s}^{(p,p')}(L,1)$, is not contained in \eqref{F}.
Using \eqref{rzero} these cases can however be obtained
from \cite[Eq. (10.2)]{BMS98} and \cite[Eq. (8.68)]{BMS98} as follows.
\begin{theorem}\label{thmF}
For $1\leq p<p'<2p$ such that $\gcd(p,p')=1$ and 
$s=l_{\sigma+1}$ a Takahashi length,
\begin{align}\label{F0s}
X_{0,s}^{(p,p')}(L,1)
&=q^{\frac{L}{2}+\Delta_s}f(L\ve_1+\vu_0,\vu_\sigma) \\ \label{F0ps}
X_{0,p'-s}^{(p,p')}(L,1)&=q^{\frac{L}{2}+\Delta'_s}
f(L\ve_1+\vu_0,\vu_\sigma+\vu_{d+1}).
\end{align}
\end{theorem}
As before, $\Delta_s$  and $\Delta'_s$ are determined by demanding 
that the left-hand side is $1$ for $q=0$, and $\vu_i$ is as defined
in equation~\eqref{vui}.

Fermionic forms for $p'>2p$ can be obtained from the previous two
theorems by the duality transformation \eqref{dual symmetry}
(and equation \eqref{rzero} when $r=0$, $b=1$).
Applying \eqref{qbininv}, this yields
\begin{equation*}
X_{b-r,s}^{(p'-p,p')}(L,b) = q^{\frac{1}{4}(L^2-(b-s)^2)-\Delta_{b,s}}
f(\vu_\sigma,L\ve_1+\vu_\beta)
\end{equation*}
and
\begin{align}\label{F0sdual}
X_{0,s}^{(p'-p,p')}(L,1)&=q^{\frac{1}{4}(L^2-s^2+1)-\Delta_s}
f(\vu_\sigma,L\ve_1+\vu_0) \\ \label{F0psdual}
X_{0,p'-s}^{(p'-p,p')}(L,1)&=q^{\frac{1}{4}(L^2-(p'-s)^2+1)-\Delta'_s}
f(\vu_\sigma+\vu_{d+1},L\ve_1+\vu_0).
\end{align}

\section{Fermionic representations of A$^{(1)}_1$ string functions
and parafermion characters}\label{secSC}

Our two main results obtained so far can be summarized as follows:
\begin{enumerate}
\item
The conjugate Bailey pairs 
$(\gamma,\delta)$ of Corollary~\ref{cor CBP} where $\gamma$ is a sequence
of (generalized) A$^{(1)}_1$ string functions and $\delta$ a sequence of 
one-dimensional configuration sums.
\item
A fermionic representation for the sequences $\delta$ as formulated
in Theorem~\ref{thmF} and equations \eqref{F0sdual} and \eqref{F0psdual}.
\end{enumerate}
As a consequence of these results 
we find fermionic or constant-sign expressions for the sequence $\gamma$
and thus for the A$^{(1)}_1$ string functions. 
Specifically, by Corollary~\ref{cor CBP} and equation \eqref{CBP} we have
\begin{equation}\label{CX}
\C_{m,\ell}^{(p,p')}(q)=\sum_{r=0}^{\infty}
\frac{X_{0,\ell+1}^{(p,p')}(2r+m,1)}{(q)_r(q)_{r+m}}
\end{equation}
and hence, using \eqref{F0s} and \eqref{F0sdual}, the following
result arises.
\begin{corollary}\label{corCF}
For $1\leq p<p'<2p$ with $\gcd(p,p')=1$ set $N=p'/p-2$, and
let $m\in\Zp$ and $\ell+1=l_{\sigma+1}$ a Takahashi length
such that $\ell+m$ is even. Then 
\begin{align}\label{CF}
\C_{m,\ell}^N(q)&=q^{\Delta_{\ell+1}+\frac{1}{2}m}
\sum_{r=0}^{\infty}
\frac{q^r f((2r+m)\ve_1+\vu_0,\vu_\sigma)}
{(q)_r(q)_{r+m}} \\
\intertext{and}
\C_{m,\ell}^{-N/(N+1)}(q)&=
q^{\frac{1}{4}(m^2-\ell(\ell+2))-\Delta_{\ell+1}}
\sum_{r=0}^{\infty}
\frac{q^{r(r+m)} f(\vu_\sigma,(2r+m)\ve_1+\vu_0)}
{(q)_r(q)_{r+m}}. \notag
\end{align}
\end{corollary}
Similarly, using  \eqref{CX}, \eqref{F0ps} and \eqref{F0psdual} we get
\begin{corollary}\label{corCF2}
For $1\leq p<p'<2p$ with $\gcd(p,p')=1$ set $N=p'/p-2$, and
let $m\in\Zp$ and $p'-\ell-1=l_{\sigma+1}$ a Takahashi length
such that $\ell+m$ is even. Then 
\begin{align*}
\C_{m,\ell}^N(q)&=q^{\Delta'_{p'-\ell-1}+\frac{1}{2}m}
\sum_{r=0}^{\infty}
\frac{q^r f((2r+m)\ve_1+\vu_0,\vu_\sigma+\vu_{d+1})}
{(q)_r(q)_{r+m}} \\
\intertext{and}
\C_{m,\ell}^{-N/(N+1)}(q)&=
q^{\frac{1}{4}(m^2-\ell(\ell+2))-\Delta'_{p'-\ell-1}} \notag \\
& \qquad \times
\sum_{r=0}^{\infty}
\frac{q^{r(r+m)} f(\vu_\sigma+\vu_{d+1},(2r+m)\ve_1+\vu_0)}
{(q)_r(q)_{r+m}}.
\end{align*}
\end{corollary}

For most choices of $p$ and $p'$ we believe these results 
to be new, but for $p=1$, $p'\geq 3$ and for $p=2$, $p'=3$ we recover
known summation formulas.
The simplest cases are $p=1$ or $2$ and $p'=3$ when we can employ 
Schur's \cite{Schur17} polynomial analogue of the Euler identity,
$X_{1,\ell+1}^{(2,3)}(L)=1$, so that by \eqref{dual symmetry} 
and \eqref{rzero}
\begin{align*}
X_{0,\ell+1}^{(1,3)}(L,1)&=q^{\frac{1}{4}(L^2-\ell^2)} \\
X_{0,\ell+1}^{(2,3)}(L,1)&=q^{\frac{1}{2}(L-\ell)}.
\end{align*}
Considering $(p,p')=(1,3)$ we find from Corollary~\ref{cor CBP}
that $\delta_L=a^L q^{L^2+(\eta^2-\ell^2)/4}$, which we recognize as
Bailey's original sequence $\delta$ of equation \eqref{gdinf}
up to an irrelevant factor $q^{(\eta^2-\ell^2)/4}$.
Hence $\gamma_L=a^L q^{L^2+(\eta^2-\ell^2)/4}/(aq)_{\infty}$ and
\begin{equation}\label{level1sf}
\C_{m,\ell}^1(q)=
\frac{q^{\frac{1}{4}(m^2-\ell^2)}}{(q)_{\infty}},
\end{equation}
which is the well-known form of the level-$1$ string 
function~\cite[Sec. 4.6, Ex. 3]{KP84}.
Next let $(p,p')=(2,3)$. Then Schur's polynomial identity implies
$\delta_L=q^{L+(\eta-\ell)/2}$ which corresponds to the specialization
$r=q$ in the sequence $\delta$ of Bressoud and Singh given in 
equation~\eqref{BS}. Accordingly, we find that
the string function at level $-1/2$ can be represented as
\begin{equation*}
\C_{m,\ell}^{-1/2}(q)=
\frac{q^{\frac{1}{2}(m-\ell)}}{(q)_{\infty}^2}
\sum_{i\in\Zp}(-1)^i q^{\frac{1}{2}i(i+2m+1)}.
\end{equation*}
A constant-sign expression can be obtained from \eqref{CF},
\begin{equation*}
\C_{m,\ell}^{-1/2}(q)
=q^{\frac{1}{2}(m-\ell)}\sum_{r=0}^{\infty}\frac{q^r}{(q)_r(q)_{r+m}}.
\end{equation*}
Using Heine's $_2\phi_1$ transformation formula \cite[Eq. (III.3)]{GR90}
\begin{equation}\label{Heine}
{_2\phi_1}\Bigl[\genfrac{}{}{0pt}{}{a,b}{c};q,z\Bigr]
=\frac{(abz/c)_{\infty}}{(z)_{\infty}}
{_2\phi_1}\Bigl[\genfrac{}{}{0pt}{}{c/a,c/b}{c};q,\frac{abz}{c}\Bigr],
\end{equation}
with $a=b=0$, $c=q^{m+1}$ and $z=q$, this can be transformed into
\begin{equation*}
\C_{m,\ell}^{-1/2}(q)
=\frac{q^{\frac{1}{2}(m-\ell)}}{(q)_{\infty}}\sum_{r=0}^{\infty}
\frac{q^{r(r+m+1)}}{(q)_r(q)_{r+m}}
\end{equation*}
which has an explicit factor $1/(q)_{\infty}$ and hence also provides
a fermionic expression for the parafermion characters
$e_{m,\ell}^{-1/2}(q)$ of equation \eqref{edef}.

By far the most involved of the known cases is $(p,p')=(1,p')$ 
for arbitrary $p'\geq 3$.
Then $N=p'-2\in\N$, $\ell\in\Z_{N+1}$, and, according to Lepowsky and
Primc~\cite{LP85},
\begin{equation}\label{LP}
\C_{m,\ell}^N(q)=\frac{q^{\frac{m^2-\ell^2}{4N}}}{(q)_{\infty}}
\sum_{\substack{\vn\in\Zp^{N-1}\\
\frac{m+\ell}{2N}+(C^{-1}\vn)_1\in\Z}}
\frac{q^{\vn C^{-1}(\vn-\ve_{\ell})}}{(q)_{\vn}},
\end{equation}
where $C$ is the A$_{N-1}$ Cartan matrix, $\ve_i$ is the $i$th 
standard unit vector in $\Z^{N-1}$ ($\ve_0=\ve_N=0$) and
$(q)_{\vn}=\prod_{j=1}^{N-1}(q)_{n_j}$.
{}From the fermionic representations
\eqref{F0sdual} and \eqref{F0psdual} for the configuration sum we also have
\begin{equation}
X^{(1,N+2)}_{0,\ell+1}(L,1)
=q^{\frac{L^2-\ell^2}{4N}}
\sum_{\substack{\vn\in\Zp^{N-1}\\
\frac{L+\ell}{2N}+(C^{-1}\vn)_1\in\Z}}
q^{\vn C^{-1}(\vn-\ve_{\ell})}
\qbin{\vm+\vn}{\vn},
\end{equation}
with $\vm+\vn=\frac{1}{2}(L\ve_1+\ve_{\ell}+\I\vm)$, and
\begin{equation}\label{X2}
X^{(1,N+2)}_{0,\ell+1}(L,1)
=q^{\frac{L^2-\ell^2}{4N}}
\sum_{\substack{\vn\in\Zp^{N-1}\\
\frac{L-\ell}{2N}+(C^{-1}\vn)_1\in\Z}}
q^{\vn C^{-1}(\vn-\ve_{N-\ell})}
\qbin{\vm+\vn}{\vn},
\end{equation}
with $\vm+\vn=\frac{1}{2}(L\ve_1+\ve_{N-\ell}+\I\vm)$.
Here $\I$ is the incidence matrix of the A$_{N-1}$ Dynkin diagram.

Inserting \eqref{LP}--\eqref{X2} into \eqref{CBP2} we obtain two sequences
of conjugate Bailey pairs. Using the symmetry
$\C_{m,\ell}^N(q)=q^{(m-\ell)/2}\C_{m-N,N-\ell}^N(q)$
these two sequences may be succinctly expressed as follows.
\begin{theorem}
For $N\geq 1$, $\sigma\in\Z_2$, $\eta\in\Zp$ and 
$\ell\in\Z_{N+1}$ such that $\ell+\eta+\sigma N$ is even,
the following pair of sequences $(\gamma,\delta)$
forms a conjugate Bailey pair relative to $a=q^{\eta}$:
\begin{align*}
\gamma_L&=\frac{a^{L/N}q^{L^2/N}}{(aq)_{\infty}}
\sum_{\substack{\vn\in\Zp^{N-1}\\
\frac{2L+\eta+\ell}{2N}+(C^{-1}\vn)_1\in\Z+\frac{\sigma}{2}}}
\frac{q^{\vn C^{-1}(\vn-\ve_{\ell})}}{(q)_{\vn}} \\
\delta_L&=a^{L/N}q^{L^2/N}
\sum_{\substack{\vn\in\Zp^{N-1}\\
\frac{2L+\eta+\ell}{2N}+(C^{-1}\vn)_1\in\Z+\frac{\sigma}{2}}}
q^{\vn C^{-1}(\vn-\ve_{\ell})}
\qbin{\vm+\vn}{\vn},
\end{align*}
with $\vm+\vn=\frac{1}{2}((2L+\eta)\ve_1+\ve_{\ell}+\I\vm)$.
\end{theorem}
These are the ``higher-level'' conjugate Bailey pairs 
of \cite[Lemma 3]{SW97} and \cite[Cor. 2.1]{SW98} (with the parameter 
$M$ therein sent to infinity and with the partition $\lambda$ therein
having a single part).

To conclude this section we give some examples of \eqref{CX} that are new.
When we take $(p,p')=(2,5)$ we can express the string functions at 
level $1/2$ in terms of polynomials introduced by 
Schur~\cite{Schur17} in his famous paper on the Rogers--Ramanujan identities.
To be specific, from~\eqref{F0sdual} we infer the following polynomial
analogues of the Rogers--Ramanujan identities
\begin{align*}
X^{(2,5)}_{0,1}(2L,1)&=q^L\Bigl(1+
\sum_{n=1}^{L-1}q^{n(n+1)}\qbins{2L-2-n}{n}\Bigr) \\
X^{(2,5)}_{0,2}(2L+1,1)&=q^L
\sum_{n=0}^L q^{n^2}\qbins{2L-n}{n} \\
X^{(2,5)}_{0,3}(2L,1)&=q^{L-1}
\sum_{n=0}^{L-1}q^{n^2}\qbins{2L-1-n}{n}\\
X^{(2,5)}_{0,4}(2L+1,1)&=q^{L-1}
\sum_{n=0}^{L-1}q^{n(n+1)}\qbins{2L-1-n}{n}.
\end{align*}
We remark that the above results may also be derived using
related polynomial identities for 
$X^{(2,5)}_{1,1}(2L,3)$,
$X^{(2,5)}_{1,1}(2L+1,2)$,
$X^{(2,5)}_{1,3}(2L,3)$ and
$X^{(2,5)}_{1,3}(2L+1,2)$,
due to Andrews~\cite{Andrews70}.
Substituting the above four identities into \eqref{CX} gives
fermionic representation for the string functions at level $1/2$.
Fermionic forms for the corresponding parafermion characters
$e_{m,\ell}^{1/2}$ can be obtained by
pulling out an explicit factor $1/(q)_{\infty}$.
\begin{proposition}
For $m\geq 0$ the level $1/2$ string functions can be expressed as
\begin{align*}
\C_{2m,0}^{1/2}(q)&=\frac{q^m}{(q)_{\infty}}
\sum_{r=0}^{\infty}\frac{q^r}{(q)_r}\Bigl(
1+\sum_{n=1}^{m+\lfloor (r-2)/2 \rfloor} q^{n(n+1)}\qbin{r+2m-n-2}{n}\Bigr)
\\[2mm]
\C_{2m+1,1}^{1/2}(q)&=\frac{q^m}{(q)_{\infty}}
\sum_{r=0}^{\infty}\frac{q^r}{(q)_r}
\sum_{n=0}^{m+\lfloor r/2 \rfloor} q^{n^2}\qbin{r+2m-n}{n} \\[2mm]
\C_{2m,2}^{1/2}(q)&=\frac{q^{m-1}}{(q)_{\infty}}
\sum_{r=0}^{\infty}\frac{q^r}{(q)_r}
\sum_{n=0}^{m+\lfloor (r-1)/2 \rfloor} q^{n^2}\qbin{r+2m-n-1}{n} \\[2mm]
\C_{2m+1,3}^{1/2}(q)&=\frac{q^{m-1}}{(q)_{\infty}}
\sum_{r=0}^{\infty}\frac{q^r}{(q)_r}
\sum_{n=0}^{m+\lfloor (r-1)/2 \rfloor} q^{n(n+1)}\qbin{r+2m-n-1}{n}.
\end{align*}
\end{proposition}
\begin{proof}
We only present the proof of the second identity. The other three identities
can be proven in a similar fashion.
(The second rather than the first identity is chosen because
all equations are more compact in this case.)
We start with 
\begin{equation*}
\C_{2m+1,1}^{1/2}(q)=q^m\sum_{r=0}^{\infty}\sum_{n=0}^{r+m} 
\frac{q^{r+n^2}}{(q)_r(q)_{r+2m+1}}\qbin{2r+2m-n}{n}
\end{equation*}
and interchange the sums over $r$ and $n$ and shift $r\to r+n-m$. 
Then we again swap the order of summation yielding
\begin{equation}\label{S1S2}
\C_{2m+1,1}^{1/2}(q)=
\Bigl(\sum_{r=m}^{\infty}\sum_{n=0}^{\infty}+
\sum_{r=0}^{m-1}\sum_{n=m-r}^{\infty}\Bigr)
\frac{q^{r+n(n+1)}}{(q)_{r+n-m}(q)_{r+n+m+1}}\qbin{n+2r}{n}.
\end{equation}
Now consider the first double sum denoted by $S_1$ and write this as
\begin{equation*}
S_1=
\sum_{r=m}^{\infty}\frac{q^r}{(q)_{r-m}(q)_{r+m+1}}
\sum_{n=0}^{\infty} \frac{q^{n(n+1)}(q^{2r+1})_n}
{(q)_n(q^{r-m+1})_n(q^{r+m+2})_n}.
\end{equation*}
Using the $q$-Kummer--Thomae--Whipple formula \cite[(III.9)]{GR90}
\begin{equation}\label{qKTW}
{_3\phi_2}\Bigl[\genfrac{}{}{0pt}{}{a,b,c}{d,e};q,\frac{de}{abc}\Bigr]
=\frac{(e/a,de/bc)_{\infty}}{(e,de/abc)_{\infty}}
{_3\phi_2}\Bigl[\genfrac{}{}{0pt}{}{a,d/b,d/c}
{d,de/bc};q,\frac{e}{a}\Bigr],
\end{equation}
with $a,b\to\infty$, $c=q^{2r+1}$, $d=q^{r-m+1}$ and $e=q^{r+m+2}$ this can
be put in the form
\begin{equation*}
S_1=\frac{1}{(q)_{\infty}}
\sum_{r=m}^{\infty}
\sum_{n=0}^{r+m} \frac{q^{r+n(n+1)}}{(q)_{r+n-m}}\qbin{r+m}{n}.
\end{equation*}
Once more the order of summation is reversed, then $r$ is replaced by $r+m-n$
and the summation order is again changed. Thus,
\begin{equation*}
S_1=\frac{q^m}{(q)_{\infty}}
\sum_{r=0}^{\infty}
\frac{q^r}{(q)_r}
\sum_{n=0}^{\min\{r,m+\lfloor r/2\rfloor\}}
q^{n^2}\qbin{r+2m-n}{n}.
\end{equation*}
Next we deal with $S_2$, given by the second  double sum in \eqref{S1S2}.
Shifting $n\to n+m-r$ gives
\begin{equation*}
S_2=\frac{q^m}{(q)_{2m+1}}\sum_{r=0}^{m-1}
q^{(m-r)^2}\qbin{r+m}{2r}
\sum_{n=0}^{\infty}
\frac{q^{n(n+2m-2r+1)}(q^{m+r+1})_n}{(q)_n(q^{m-r+1})_n(q^{2m+2})_n}.
\end{equation*}
By equation \eqref{qKTW} with $a,b\to\infty$, $c=q^{m+r+1}$, $d=q^{m-r+1}$
and $e=q^{2m+2}$ this is equal to
\begin{equation*}
S_2=\frac{1}{(q)_{\infty}}\sum_{r=0}^{m-1}\sum_{n=0}^{2r}
\frac{q^{r+(n+m-r)(n+m-r+1)}}{(q)_n}\qbin{r+m}{n+m-r}.
\end{equation*}
By an interchange of sums followed by the successive transformations 
$r\to n+m-r$ and $r\leftrightarrow n$ this becomes
\begin{equation*}
S_2=\frac{q^m}{(q)_{\infty}}\sum_{r=0}^{2m-2}\frac{q^r}{(q)_r}
\sum_{n=r+1}^{m+\lfloor r/2\rfloor}q^{n^2}\qbin{r+2m-n}{n}.
\end{equation*}
Computing $S_1+S_2$ results in the claim of the proposition.
\end{proof}

In our last example we take $(p,p')=(3,4)$.
The one-dimensional configuration sums for this case correspond to those 
of the celebrated Ising model of statistical mechanics, and the fermionic 
representations of the previous section can be simplified using the
$q$-binomial theorem \eqref{qbinthm} or the $q$-Chu--Vandermonde sum
\eqref{qCV}. Specifically we have the polynomial identities
\begin{equation}\label{p34}
X_{0,1}^{(3,4)}(2L,1)\pm q^{3/2}X_{0,3}^{(3,4)}(2L,1)=q^L(\mp q^{1/2})_L
\end{equation}
and
\begin{equation}\label{p34'}
X_{0,2}^{(3,4)}(2L+1,1)=q^L(-q)_L.
\end{equation}
Substitution into \eqref{CX} yields fermionic forms for the
string functions at level $-2/3$. The next proposition states alternative
expressions for these string functions which by \eqref{edef}
also imply fermionic forms for the corresponding parafermion characters. 
\begin{proposition}
For $m\geq 0$ the level $-2/3$ string functions satisfy the identities
\begin{align*}
\C_{2m,0}^{-2/3}(q)
&=\frac{q^m}{2(q)_{\infty}}
\sum_{r=0}^{\infty}\frac{q^{r^2/2+(m+1)r}}{(q)_r(q)_{r+2m}}
\{(-q^{1/2})_{r+m}+(-1)^r (q^{1/2})_{r+m}\}\\
q^{3/2}\C_{2m,2}^{-2/3}(q)
&=\frac{q^m}{2(q)_{\infty}}
\sum_{r=0}^{\infty}\frac{q^{r^2/2+(m+1)r}}{(q)_r(q)_{r+2m}}
\{(-q^{1/2})_{r+m}-(-1)^r (q^{1/2})_{r+m}\}\\
\C_{2m+1,1}^{-2/3}(q)
&=\frac{q^m}{(q)_{\infty}}
\sum_{r=0}^{\infty}\frac{q^{\binom{r}{2}+(m+2)r}(-q)_{r+m}}
{(q)_r(q)_{r+2m+1}}.
\end{align*}
\end{proposition}
\begin{proof}
Inserting the polynomial identities \eqref{p34} and \eqref{p34'}
into \eqref{CX} one can apply the $_2\phi_1$ transformation \eqref{Heine} 
with $a=0$, $b=\mp q^{m+1/2}$, $c=q^{2m+1}$, $z=q$, and $a=0$, $b=-q^{m+1}$, 
$c=q^{2m+2}$, $z=q$, respectively. This yields identities for
$\C_{2m,0}^{-2/3}(q)\pm q^{3/2}\C_{2m,2}^{-2/3}(q)$ and
$\C_{2m+1,1}^{-2/3}(q)$ which immediately imply the expressions of 
the proposition.
\end{proof}
Note that one can apply \eqref{qbinthm} once again to rewrite
$$
\frac{1}{2}\{(-q^{1/2})_{r+m}\pm (-1)^r(q^{1/2})_{r+m}\}
=\sum_{n,\text{restriction}} q^{n^2/2}\qbin{r+m}{n},
$$
where the restrictions are $n\equiv r \pmod 2$ and
$n\not\equiv r \pmod 2$, respectively.

\section{A$_1^{(1)}$ branching functions}\label{sec branchingf}

Let either $N_1$ or $N_2$ be a positive integer. Then the A$_1^{(1)}$
branching functions are defined by~\cite{KW90}
\begin{equation}\label{bf def}
\chi_{\ell_1}^{N_1}(z,q)\chi_{\ell_2}^{N_2}(z,q)
=\sum_{\substack{\ell_3\in \Z_{p_3'-1} \\ \ell_1+\ell_2+\ell_3\in 2\Z}}
B_{\ell_1,\ell_2,\ell_3}^{N_1,N_2}(q)
\chi_{\ell_3}^{N_3}(z,q).
\end{equation}
Here $N_1=p'_1/p_1-2$, $N_2=p'_2/p_2-2$ and $N_3=N_1+N_2=p'_3/p_3-2$,
with $\gcd(p_i,p_i')=1$ for $i=1,2,3$.
Note that $p_3=p_1 p_2$ and $p_3'=p_1'p_2+p_2'p_1-2p_1p_2=p_2(p_1'+N_2p_1)$.
Indeed $\gcd(p_3,p_3')=1$ since either $p_1=1$ or $p_2=1$.

In the following we are going to derive an explicit expression for
the branching function following the method employed by Kac and Wakimoto
in \cite{KW90} (see also~\cite{DJKMO87,DJKMO88}). 
The essence of this approach is to expand the character
$\chi_{\ell_2}^{N_2}$ in terms of string functions and to then perform
simple manipulations using the symmetries of the string functions
to express the left-hand side of \eqref{bf def} as a linear combination
of the $\chi_{\ell_3}^{N_3}$. The difference between our derivation below
and that of Kac and Wakimoto is that we will not assume that $N_2$
is integer. Of course, since either $N_1$ or $N_2$ is (a positive) integer and
$B^{N_1,N_2}_{\ell_1,\ell_2,\ell_3}=B^{N_2,N_1}_{\ell_2,\ell_1,\ell_3}$
one can without loss of generality assume that $N_2\in\N$. Nevertheless,
dropping this assumption leads to a different representation
of the branching functions. As will be shown in the next section, this
has a natural interpretation in terms of the Bailey lemma.
Before we commence our derivation we remark that because $N_2$ is no longer
assumed to be integer we deal with string functions at (generally) 
non-integer level and hence we cannot rely on the symmetries employed
in the Kac--Wakimoto derivation.

Insert \eqref{char} for $\chi_{\ell_1}^{N_1}(z,q)$ 
and \eqref{sfdef} for $\chi_{\ell_2}^{N_2}(z,q)$ in the left-hand side
of \eqref{bf def}.
Then, using the definition \eqref{Thetadef} of $\Theta_{n,m}(z,q)$, 
one obtains
\begin{align}\label{P}
P^{N_1,N_2}_{\ell_1,\ell_2}(q)&:=
\chi_{\ell_1}^{N_1}(z,q)\chi_{\ell_2}^{N_2}(z,q)
\sum_{\sigma=\pm 1}\sigma\Theta_{\sigma,2}(z,q) \\
&\,=\sum_{\sigma=\pm 1}\sum_{j\in\Z+\sigma\frac{\ell_1+1}{2p'_1}}
\sum_{m\in 2\Z+\ell_2}
\sigma z^{-\frac{1}{2}(m+2p'_1j)}q^{\frac{m^2}{4N_2}+p_1 p_1'j^2}
C_{m,\ell_2}^{N_2}(q). \notag 
\end{align}
Now make the replacement $m\to m-2p_1'j$ followed by 
$j\to \sigma(j+\frac{\ell_1+1}{2p_1'})$. Using $C_{m,\ell}^N=C_{-m,\ell}^N$
this gives
\begin{equation}\label{inter}
P^{N_1,N_2}_{\ell_1,\ell_2}(q)=
q^{\frac{(\ell_1+1)^2}{4(N_1+2)}}
\sum_{m\in 2\Z+\ell_1+\ell_2+1}
z^{-\frac{1}{2}m}
q^{\frac{1}{4N_2}(m-\ell_1-1)^2}
b_{\ell_1+1,m,\ell_2}^{p_1',p_1'+N_2p_1,N_2}(q),
\end{equation}
where we have introduced the function
\begin{multline*}
b_{r,s,\ell}^{P,P',N}(q) \\=
\sum_{j\in\Z}\Bigl\{q^{\frac{j}{N}(PP'j+P'r-Ps)}C_{2Pj+r-s,\ell}^N(q)
-q^{\frac{1}{N}(Pj+r)(P'j+s)}C_{2Pj+r+s,\ell}^N(q)\Bigr\}.
\end{multline*}
Note that the initial assumption that either $N_1$ or $N_2$ is a positive
integer means that we are only concerned with 
$b_{r,s,\ell}^{P,P',N}(q)$ with either $(P'-P)/N=1$ or $N\in\N$.
This is crucial in the following lemma needed to 
rewrite the expression for $P^{N_1,N_2}_{\ell_1,\ell_2}(q)$.

\begin{lemma}\label{lem periodic}
Let $P\in\N$ and $N,P'\in\Q$ such that
$N=p'/p-2$ with $\gcd(p,p')=1$ and $(P'-P)/N\in\Zp$.
When $(P'-P)/N=1$ or $N\in\N$ the following periodicity holds:
\begin{equation}\label{periodicity}
b_{r,s+2pP',\ell}^{P,P',N}(q)=
q^{-\frac{p}{N}(pPP'-P'r+Ps)}b_{r,s,\ell}^{P,P',N}(q).
\end{equation}
\end{lemma}

\begin{proof}
After inserting the definition of $b_{r,s,\ell}^{P,P',N}$ in the above 
equation make the variable changes $j\to j+p$ in the first term and
$j\to j-p$ in the second term of the left-hand side.
Then, by the symmetry \eqref{sym1}, equation \eqref{periodicity} can
be rewritten as
\begin{multline}\label{Ceq}
\sum_{j\in\Z}\Bigl\{
q^{\frac{j}{N}(PP'j+P'r-Ps)}C_{2Pj+r-s-2pkN,\ell}^N(q)
-q^{\frac{1}{N}(Pj-r)(P'j-s)}C_{2Pj-r-s-2pkN,\ell}^N(q)\Bigr\} \\
\sum_{j\in\Z}\Bigl\{
q^{\frac{j}{N}(PP'j+P'r-Ps)}C_{2Pj+r-s,\ell}^N(q)
-q^{\frac{1}{N}(Pj-r)(P'j-s)}C_{2Pj-r-s,\ell}^N(q)\Bigr\}
\end{multline}
where $k=(P'-P)/N\in\Zp$.
When $N\in\N$ this follows directly from the symmetries
\eqref{sym1} and \eqref{sym2}, and in the remainder we assume that 
$N\in\Q$ and $k=1$.
The complication is now that we no longer have $C^N_{m,\ell}=C^N_{m-2N,\ell}$.
In view of this let us first investigate the origin of this difficulty.
Consider the expression \eqref{sffinal} of the A$_1^{(1)}$ string functions.
The summand has two different terms corresponding to the two terms within
the curly braces. In the first term make the variable change
$j\to j-1$, $i\to i+2p$ and in the second term make the change
$j\to j+1$, $j\to i+2p$.
The result of these changes is exactly the same expression as before
except that $m$ has been replaced by $m-2pN$ and that the sum over $i$
now runs over all integers greater than $-2p$.
We may therefore conclude that
\begin{equation}\label{CCb}
C_{m,\ell}^N(q)=C_{m-2pN,\ell}^N(q)+\bar{C}_{m-2pN,\ell}^{N}(q),
\end{equation}
where
\begin{multline*}
\bar{C}_{m,\ell}^N(q)=
\frac{q^{\frac{(\ell+1)^2}{4(N+2)}-\frac{m^2}{4N}}}
{\eta^3(\tau)} \\ \times \sum_{i=1}^{2p-1}\sum_{j\in\Z}
(-1)^i q^{\frac{1}{2}i(i-m)+pj(p'j+\ell+1)}
\Bigl\{q^{-\frac{1}{2}i(2p'j+\ell+1)}-q^{\frac{1}{2}i(2p'j+\ell+1)}\Bigr\}.
\end{multline*}
By a shift $j\to j-1$ in the second term of the summand this becomes
\begin{multline*}
\bar{C}_{m,\ell}^N(q)=
\frac{q^{\frac{(\ell+1)^2}{4(N+2)}-\frac{m^2}{4N}}}
{\eta^3(\tau)} \\ \times \sum_{i=1}^{2p-1}\sum_{j\in\Z}
(-1)^i q^{\frac{1}{2}i(i-m)+pj(p'j+\ell+1)-\frac{1}{2}i(2p'j+\ell+1)}
\Bigl\{1-q^{\frac{1}{2}(i-p)(2p'j-p'+\ell+1)}\Bigr\},
\end{multline*}
which shows that the $i=p$ term in the summand vanishes and hence that
$\bar{C}_{m,\ell}^N(q)=0$ for $N$ integer.

Inserting \eqref{CCb} into equation \eqref{Ceq} with $k=1$
we are done with the lemma if we prove that
\begin{multline*}
\sum_{j\in\Z}\Bigl\{q^{\frac{j}{N}(PP'j+P'r-Ps)}
\bar{C}_{2Pj+r-s-2pN,\ell}^N(q) \\
-q^{\frac{1}{N}(Pj-r)(P'j-s)}
\bar{C}_{2Pj-r-s-2pN,\ell}^N(q)\Bigr\}=0.
\end{multline*}
Using the explicit form for $\bar{C}_{m,\ell}^N(q)$, this is equivalent 
to showing that
\begin{multline*}
\sum_{i=1}^{2p-1}\sum_{j\in\Z}
(-1)^i q^{\frac{1}{2}i(i-r+s)+jp(p'j+\ell+1)
-\frac{1}{2}i(2p'j+\ell+1)} 
\Bigl\{1-q^{\frac{1}{2}(i-p)(2p'j-p'+\ell+1)}\Bigr\} \\
\times q^{p((i-p)N+r-s)}
\sum_{\mu\in\Z}\Bigl\{
q^{\mu(\mu P+r-P(i-2p))}
-q^{(\mu-i+2p)(\mu P-r)}\Bigr\}=0.
\end{multline*}
After the shift $\mu\to i-2p-\mu$ in the second term in the sum over $\mu$
we are done.
\end{proof}

{}From \eqref{sym1} it follows that $b_{r,s,\ell}^{P,P,N}(q)=
-q^{\frac{rs}{N}}b_{r,-s,\ell}^{P,P',N}(q)$ so that in combination
with Lemma \ref{lem periodic} 
\begin{equation}\label{periodicity1}
b_{r,2pP'-s,\ell}^{P,P',N}(q)=
 -q^{-\frac{1}{N}(pP-r)(pP'-s)}b_{r,s,\ell}^{P,P',N}(q).
\end{equation}
In view of \eqref{periodicity} and \eqref{periodicity1}, it becomes 
natural to dissect the sum over $m$ in \eqref{inter} using
\begin{equation*}
\sum_{m\in 2\Z+\ell_1+\ell_2+1} f_m=
\sum_{k\in\N}\Bigl\{
\sum_{\substack{\ell_3\in\Z_{p_3'} \\ \ell_1+\ell_2+\ell_3\in 2\Z}}
f_{2p_3'k+\ell_3+1}+
\sum_{\substack{\ell_3+1\in\Z_{p_3'} \\ \ell_1+\ell_2+\ell_3\in 2\Z}}
f_{2p_3'k-\ell_3-1}\Bigr\}.
\end{equation*}
Observing that $b_{r,0,\ell}^{P,P',N}=0$ and, by \eqref{periodicity1},
also $b_{r,pP',\ell}^{P,P',N}=0$, equation \eqref{inter}
can then be written as
\begin{multline*}
P^{N_1,N_2}_{\ell_1,\ell_2}(q)=
q^{\frac{(\ell_1+1)^2}{4(N_1+2)}+\frac{(\ell_3-\ell_1)^2}{4N_2}}
\sum_{\substack{\ell_3\in\Z_{p'_3-1}\\\ell_1+\ell_2+\ell_3\in2\Z}}
b_{\ell_1+1,\ell_3+1,\ell_2}^{p_1',p_1'+N_2 p_1,N_2}\\
\quad\times \sum_{k\in\Z} \Bigl\{
 z^{-\frac{1}{2}(2p'_3k+\ell_3+1)} q^{p_3k(p'_3k+\ell_3+1)} 
-z^{-\frac{1}{2}(2p'_3k-\ell_3-1)} q^{p_3k(p'_3k-\ell_3-1)}\Bigr\}\\
=q^{\frac{((p_1'+N_2 p_1)(\ell_1+1)-p_1'(\ell_3+1))^2}{4N_2p_1'(p_1'+N_2 p_1)}}
\sum_{\substack{\ell_3\in\Z_{p'_3-1}\\\ell_1+\ell_2+\ell_3\in2\Z}}
b_{\ell_1+1,\ell_3+1,\ell_2}^{p_1',p_1'+N_2 p_1,N_2}
\sum_{\sigma=\pm 1}\sigma \Theta_{\sigma(\ell_3+1),p'_3}(z,q^{p_3}).
\end{multline*}
Comparing with \eqref{bf def} and \eqref{P} we can read off the
branching functions.
\begin{theorem}
For $N_1=p_1'/p_1-2$ and $N_2=p_2'/p_2-2$ with 
$\gcd(p_1,p_1')=\gcd(p_2,p_2')=1$, such that $p_1=1$ or $p_2=1$ we have
\begin{multline*}
B_{r-1,\ell,s-1}^{N_1,N_2}(q)=B_{\ell,r-1,s-1}^{N_2,N_1}(q)=
q^{\frac{(P'r-Ps)^2}{4N_2PP'}}\\[2mm]
\times \sum_{j\in\Z}\Bigl\{
q^{\frac{j}{N_2}(PP'j+P'r-Ps)}C_{2Pj+r-s,\ell}^{N_2}(q)-
q^{\frac{1}{N_2}(Pj+r)(P'j+s)}C_{2Pj+r+s,\ell}^{N_2}(q)\Bigr\}.
\end{multline*}
Here $P=p'_1$, $P'=p_1'+N_2 p_1$, $r\in\N_{P-1}$,
$\ell+1\in\N_{p'_2-1}$ and $s\in\N_{p_2 P'-1}$.
\end{theorem}
When $N_2\in\N$ this is Theorem~3.1 of \cite{KW90} for $\text{X}^{(r)}_N=
\text{A}_1^{(1)}$.

For comparison with later expressions it will be convenient to normalize
the branching functions and to express them in terms of the modified
string functions. Hence we introduce
\begin{equation}\label{BFmod}
\B_{r-1,\ell,s-1}^{N_1,N_2}(q)=
\sum_{j\in\Z}q^{p_1j(p_1'j+r)}
\Bigl\{\C_{2p_1'j+r-s,\ell}^{N_2}(q)-\C_{2p_1'j+r+s,\ell}^{N_2}(q)\Bigr\},
\end{equation}
where $$B_{r-1,\ell,s-1}^{N_1,N_2}(q)=
q^{\frac{(P'r-Ps)^2}{4N_2PP'}+\frac{(\ell+1)^2}{4(N_2+2)}
-\frac{(r-s)^2}{4N_2}-\frac{1}{8}}
\B_{r-1,\ell,s-1}^{N_1,N_2}(q).$$

\section{Bose-Fermi identities}\label{sec BF}
In Section~\ref{secSC} we have applied Corollary~\ref{cor CBP} to derive
fermionic representations for the A$_1^{(1)}$ string functions, but so far
we have not yet employed the result of Corollary~\ref{cor CBP} in the
context of the Bailey lemma. This is what we will do next.
To simplify the notation we abbreviate the polynomial
identities \eqref{F}--\eqref{F0s} as
\begin{equation}\label{XF}
X_{r,s}^{(p,p')}(L,b)=F_{r,s}^{(p,p')}(L,b).
\end{equation}
{}From these, Bailey pairs relative to $q^{|b-s|}$ can be 
extracted~\cite{Andrews84,FQ96}.
Together with the conjugate Bailey pairs of Corollary~\ref{cor CBP}
these Bailey pairs (given by \cite[Eq. (3.6)]{BMSW97})
may be substituted into Bailey's equation \eqref{abcd}.
Omitting the details we find the following theorem.
\begin{theorem}\label{thm BFid}
For $i=1,2$, let $1\leq p_i<p_i'<2p_i$ such that $\gcd(p_i,p_i')=1$ and
set $N_i=p_i'/p_i-2$.
Let $b$ and $s$ be Takahashi lengths with respect to the continued fraction 
decomposition of $p_1/(p_1'-p_1)$ and let $r=\bar{b}$.
Let $\ell+1$ be a Takahashi length with respect to the continued fraction 
decomposition of $p_2/(p_2'-p_2)$.
Then, for $\eta=|b-s|$ with $\eta+\ell$ even,
\begin{multline}\label{BFid}
\sum_{j\in\Z}\Bigl\{
q^{j(p_1 p_1'j+r p_1'-s p_1)}
\C_{2p_1'j+b-s,\ell}^{N_2}(q)-
q^{(p_1 j+r)(p_1'j+s)}
\C_{2p_1'j+b+s,\ell}^{N_2}(q)\Bigr\} \\
=\sum_{L=0}^{\infty} F_{r,s}^{(p_1,p_1')}(2L+\eta,b)
F_{0,\ell+1}^{(p_2,p_2')}(2L+\eta,1)/(q)_{2L+\eta}.
\end{multline}
\end{theorem}
Many similar theorems can be derived.
For example, we could have iterated the Bailey
pair implied by \eqref{XF} (see \cite[Eq. (3.8)]{BMSW97})
before substituting it into \eqref{abcd}.
Alternatively one can derive identities for $N_1>0$, $N_2<0$,
or $N_1<0$, $N_2>0$ or $N_1,N_2>0$.

In general we have not been able to identify the left-hand side of
\eqref{BFid}, but when either $N_1$ or $N_2$ is a positive integer 
one can recognize the left-hand side of the above identities as 
A$_1^{(1)}$ branching function.
First assume $N_2$ is integer and $r$ is even. 
Using the symmetries $\C_{m-2N_2,\ell}^{N_2}(q)
=q^{N_2-m}\C_{m,\ell}^{N_2}(q)$ and
$\C_{m,\ell}^{N_2}(q)=\C_{-m,\ell}^{N_2}(q)$, 
the left-hand side of \eqref{BFid} becomes
\begin{equation*}
q^{\frac{1}{4}r(2s-2b-N_2 r)}
\B_{s-1,\ell,b+N_2 r-1}^{N_1,N_2}(q).
\end{equation*}
For $N_2$ integer and $r$ odd we can use
$\C_{m-N_2,\ell}^{N_2}(q)
=q^{(N_2-m-\ell)/2}\C_{m,N_2-\ell}^{N_2}(q)$ and
$\C_{m,\ell}^{N_2}(q)=\C_{-m,\ell}^{N_2}(q)$
to rewrite the left-hand side of \eqref{BFid} as
\begin{equation*}
q^{\frac{1}{4}r(2s-2b-N_2 r)+\frac{1}{4}(N_2-2\ell)}
\B_{s-1,N_2-\ell,b+N_2 r-1}^{N_1,N_2}(q).
\end{equation*}
Finally, for $N_1$ integer we must have $r=0$, $b=1$ or $r=1$, 
$b-1\in\N_{p_1'-2}$
and the left-hand side of \eqref{BFid} can be simplified to
\begin{equation*}
\B_{s-1,\ell,0}^{N_1,N_2}(q) \qquad \text{and} \qquad
\B_{p_1'-s-1,\ell,p_1'-b-1}^{N_1,N_2}(q),
\end{equation*}
respectively.

Given the above results let us connect to the discussion in 
Sections~\ref{sec BL} and \ref{sec branchingf} on the duality between 
Bailey and conjugate Bailey pairs and on the symmetry of the branching
functions.
If $r=0$ and $b=1$ the right-hand side of \eqref{BFid} is symmetric under
the simultaneous interchange $N_1 \leftrightarrow N_2$ and 
$s \leftrightarrow \ell+1$.
In terms of Bailey and conjugate Bailey pairs this
corresponds to the transformation
\begin{equation*}
(\beta^{(N_1)},\delta^{(N_2)})\leftrightarrow
(\bar{\beta}^{(N_2)},\bar{\delta}^{(N_1)})
\end{equation*}
with
\begin{align*}
\bar{\beta}_L^{(N_2)}&=\delta_L^{(N_2)}/(q)_{2L+\eta} \\
\bar{\delta}_L^{(N_1)}&=\beta_L^{(N_1)}(q)_{2L+\eta},
\end{align*}
where $\beta_L^{(N_1)}=F_{0,s}^{(p_1,p_1')}(2L+\eta,1)/(q)_{2L+\eta}$
and $\delta_L^{(N_2)}=F_{0,\ell+1}^{(p_2,p_2')}(2L+\eta,1)$.
This result is to be compared with \eqref{duality}.

Similarly, using \eqref{rzero}, the right-hand side of \eqref{BFid}
is symmetric under the interchange $N_1 \leftrightarrow N_2$ and
$s \leftrightarrow \ell+1$ if $r=1$ and $b=1$,
which corresponds to the transformation
\begin{align*}
\bar{\beta}_L^{(N_2)}&=q^{-L-(\eta-\ell)/2}\delta_L^{(N_2)}/(q)_{2L+\eta} \\
\bar{\delta}_L^{(N_1)}&=q^{L+(\eta-\ell)/2}\beta_L^{(N_1)}(q)_{2L+\eta},
\end{align*}
where $\beta_L^{(N_1)}=q^{-L-(\eta-\ell)/2}
F_{0,s}^{(p_1,p_1')}(2L+\eta,1)/(q)_{2L+\eta}$
and $\delta_L^{(N_2)}=F_{0,\ell+1}^{(p_2,p_2')}(2L+\eta,1)$.

Carrying out the corresponding transformations on $\alpha$ and
$\gamma$ yields another expression for the left-hand side of~\eqref{BFid} 
which involves the modified string functions at level $N_1$.
When either $N_1$ or $N_2$ is a positive integer we recognize
the resulting identities as the special cases $\ell_3=0$ or $\ell_3=N_1$
of the symmetry
$B_{\ell_1,\ell_2,\ell_3}^{N_1,N_2}=B_{\ell_2,\ell_1,\ell_3}^{N_2,N_1}$,
as expected.

Finally we present some explicit identities that follow by application of 
Bailey's lemma and the conjugate Bailey pairs of Corollary~\ref{cor CBP}.
In Refs.~\cite[Eqs. (2.12), (2.13)]{Andrews84} and 
\cite[Eqs. (3.47), (3.48)]{Andrews85} 
one can find the following generalization of \eqref{initial},
\begin{equation}\label{BPin}
\alpha_L=\frac{(1-aq^{2L})(a)_L(-1)^L q^{\binom{L}{2}}}{(1-a)(q)_L}
\qquad \text{and} \qquad \beta_L=\delta_{L,0}.
\end{equation}
Inserting this and \eqref{CBP2} into equation \eqref{abcd} and 
performing some series manipulations gives a generalized Euler identity
for the modified string functions.
\begin{proposition}\label{propEuler}
For $1\leq p<p'$, $\ell\in\Z_{p'-1}$, $\eta\in\Z_{p'}$ such that
$\ell+\eta$ is even,
\begin{equation*}
\sum_{L=-\infty}^{\infty} (-1)^L q^{\binom{L}{2}}\C_{2L+\eta,\ell}^{(p,p')}(q)
=\delta_{\ell,\eta}.
\end{equation*}
\end{proposition}
Recalling \eqref{level1sf}, this is the classical Euler identity for
$(p,p')=(1,3)$.
For $p=1$ and arbitrary $p'$ this is the A$_1^{(1)}$ case of
equation (2.1.17) of Ref.~\cite{KW90}.

Before we can proof the proposition we need a technical lemma.
\begin{lemma}\label{lem tech}
If $f_m=f_{-m}$ then
\begin{multline}\label{lemt}
\sum_{L=0}^{\infty}
(1-q^{2L+\eta})(q^{L+1})_{\eta-1}(-1)^L q^{\binom{L}{2}} f_{2L+\eta} \\
= \sum_{k=0}^{\lfloor \eta/2 \rfloor}
\Bigl\{\qbin{\eta}{k}-\qbin{\eta}{k-1}\Bigr\}
\sum_{L=-\infty}^{\infty}(-1)^L q^{\binom{L}{2}} f_{2L+\eta-2k}.
\end{multline}
\end{lemma}
\begin{proof}
First observe that
\begin{equation*}
\sum_{k=0}^{\lfloor \eta/2 \rfloor}
\Bigl\{\qbin{\eta}{k}-\qbin{\eta}{k-1}\Bigr\}
\sum_{L=k-\eta+1}^{k-1}(-1)^L q^{\binom{L}{2}} f_{2L+\eta-2k}=0.
\end{equation*}
To prove this shift $L\to L+k$ in the first term in the curly braces 
and successively $k\to\eta-k+1$ and $L\to\eta-L+1$ in the second term 
in the curly braces. Using the symmetry of $f_m$ the resulting terms
can be combined to
\begin{equation*}
\sum_{L=1}^{\eta-1} f_{2L-\eta}\sum_{k=0}^{\eta}
(-1)^{k-L}q^{\binom{k-L}{2}}\qbin{\eta}{k}=
\sum_{L=1}^{\eta-1} f_{2L-\eta}(-1)^L q^{\binom{L+1}{2}}
(q^{-L})_{\eta}=0,
\end{equation*}
where the middle term follows by application of the $q$-binomial 
theorem \eqref{qbinthm} and the last term by $(q^{-a})_b=0$ for $0\leq a<b$.
With this result we can write the sum over $L$ in the right-hand side of
equation \eqref{lemt} as a sum over $L\leq k-\eta$ and $L\geq k$.
Then using the symmetry of $f_m$ the right-hand side becomes
\begin{equation*}
\begin{split}
& \sum_{L=0}^{\infty} f_{2L+\eta}
\sum_{k=0}^{\lfloor \eta/2 \rfloor}
\Bigl\{\qbin{\eta}{k}-\qbin{\eta}{k-1}\Bigr\}
\Bigl\{(-1)^{L+k}q^{\binom{L+k}{2}}+(-1)^{k-\eta-L}
q^{\binom{k-\eta-L}{2}}\Bigr\} \\
& \qquad \qquad = \sum_{L=0}^{\infty} f_{2L+\eta}
\sum_{k=0}^{\eta}(-1)^{L+k}q^{\binom{L+k}{2}}(1+q^{L+k})\qbin{\eta}{k} \\
& \qquad \qquad = \sum_{L=0}^{\infty} f_{2L+\eta}
(-1)^{L}q^{\binom{L}{2}}\bigl\{(q^L)_{\eta}+q^L (q^{L+1})_{\eta} \bigr\}.
\end{split}
\end{equation*}
Comparing with the left-hand side of \eqref{lemt} we are done since
$(a)_n+a(aq)_n=(1-a^2q^n)(aq)_{n-1}$.
\end{proof}

\begin{proof}[Proof of Proposition~\ref{propEuler}]
Inserting \eqref{BPin} and \eqref{CBP2} into equation \eqref{abcd}
gives the identity
\begin{equation*}
X_{0,\ell+1}^{(p,p')}(\eta,1)=
\sum_{L=0}^{\infty}
(1-q^{2L+\eta})(q^{L+1})_{\eta-1}(-1)^L q^{\binom{L}{2}}
\C_{2L+\eta,\ell}^{(p,p')}(q),
\end{equation*}
for $\eta+\ell$ even and $\ell+1\in\N_{p'-1}$.
Applying Lemma~\ref{lem tech} this can be simplified to 
\begin{equation*}
X_{0,\ell+1}^{(p,p')}(\eta,1)=
\sum_{k=0}^{\lfloor \eta/2 \rfloor}
\Bigl\{\qbin{\eta}{k}-\qbin{\eta}{k-1}\Bigr\}
\sum_{L=-\infty}^{\infty}(-1)^L q^{\binom{L}{2}}
\C_{2L+\eta-2k,\ell}^{(p,p')}(q).
\end{equation*}
Now observe that for $\eta\leq p'-1$ the only contribution to
$X_{0,\ell+1}^{(p,p')}(\eta,1)$ comes from the $j=0$ term in the
summand of \eqref{odcs}. Therefore,
\begin{equation*}
X_{0,\ell+1}^{(p,p')}(\eta,1)
=\qbin{\eta}{(\eta+\ell)/2}-\qbin{\eta}{(\eta-\ell-2)/2}
=\sum_{k=0}^{\lfloor \eta/2 \rfloor}
\Bigl\{\qbin{\eta}{k}-\qbin{\eta}{k-1}\Bigr\}\delta_{\eta-2k,\ell}.
\end{equation*}
By induction on $\eta$ this implies Proposition~\ref{propEuler}.
\end{proof}

Our last identity follows by a straightforward generalization of the proof
of Theorem~4.1 of Ref.~\cite{SW98}, which corresponds to $p=1$ in the 
result given below.
\begin{theorem}
For $1\leq p<p'$, $\ell\in\Z_{p'-1}$ and integers $\delta,k,i$ such that
$\delta\in\Z_2$, $k\geq 2$ and $i\in\N_k$,
\begin{multline}\label{RRppp}
\sum_{L=-\infty}^{\infty}(-1)^L q^{((2k+\delta-2)L+2k-2i+\delta)L/2}
\C_{2L,\ell}^{(p,p')}(q) \\
=\sum_{n_1,\dots,n_{k-1}\geq 0}
\frac{q^{N_2^2+\cdots+N_{k-1}^2+N_i+\cdots+N_{k-1}}
X_{0,\ell+1}^{(p,p')}(2N_1,1)}
{(q)_{n_1}\cdots(q)_{n_{k-2}}
(q^{2-\delta};q^{2-\delta})_{n_{k-1}}},
\end{multline}
where $N_j=n_j+\cdots+n_{k-1}$.
\end{theorem}
By Jacobi's triple product identity \eqref{tpi} and the fermionic expressions
for the string function and configuration sums given earlier in the paper,
the above identities can be
recognized as (i) Andrews' analytic counterpart of Gordon's partition theorem
when $(p,p')=(1,3)$ and $\delta=1$ \cite{Andrews74},
(ii) Bressoud's generalization thereof
to even moduli when $(p,p')=(1,3)$ and $\delta=0$ \cite{Bressoud80b}, 
(iii) generalizations of
the G\"ollnitz--Gordon partition identities due to Andrews and Bressoud
when $(p,p')=(1,4)$ and $\delta=1$ \cite{Andrews76b,Bressoud80c},
(iv) Rogers--Ramanujan type identities
by Bressoud when $(p,p')=(1,4)$ and $\delta=0$ \cite{Bressoud80c}.

\subsection*{Acknowledgements}
We thank Victor G.~Kac and Mark Shimozono for helpful comments.
The first author was supported by the ``Stichting Fundamenteel
Onderzoek der Materie''. The second author was
supported by a fellowship of the Royal
Netherlands Academy of Arts and Sciences
and a travel grant of the 
Netherlands Organization for Scientific Research (NWO).

\end{document}